%% file: Main.tex
\documentclass[11pt,twoside]{article}
\usepackage{fullpage}
\usepackage{psfig}
\usepackage{amsfonts}
\usepackage{amssymb}
\input{inputlatex}

\input{imsmark}

\renewcommand{\Bbb}[1]{{\mathbb #1}}
\newcommand{\closure}[1]{\overline{#1}}
\renewcommand{\remark}{\noindent{\bf Remark:\ }}
\renewcommand{\proof}{\noindent{\bf Proof:\ }}
\newcommand{\dfn}[1]{{\sffamily\itshape{#1}}}
\newcommand{\DS}[1]{\displaystyle{#1}}
\newcommand{\set}[1]{\left\{{#1}\right\}}
\newcommand{\st}{\,\,\vrule{}\,\,}
\newcommand{\union}{\cup}
\newcommand{\Union}{\bigcup}
\newcommand{\intersect}{\cap}

\newcommand{\text}[1]{\;\mbox{#1}\;}
\newcommand{\GamDisk}[1]{\Gamma_{\kern -2pt\odot}(#1)}
\renewcommand{\QED}{\hfill \mbox{$\square$} \vspace{6mm}}
\newcommand{\itemsetup}{%
     \setlength{\itemsep}{2pt plus 1pt minus 1pt}
     \setlength{\topsep}{2pt plus 1pt minus 1pt}
}
\newcommand{\Caption}[1]{{\caption{\footnotesize #1}}}

\advance\marginparwidth by -\textwidth          
\advance\marginparwidth by -1 true in           
\advance\marginparwidth by -\oddsidemargin      
\advance\marginparwidth by -\marginparsep       
\advance\marginparwidth by -.2 true in         

\newlength{\awidth}
\newcommand{\abox}[5]{%
  \parbox{\awidth}{\makebox[\awidth]{#1}\\
    {\footnotesize \makebox[\awidth]{\tt #5}\\
    \makebox[\awidth]{#2}\\[-3pt] \makebox[\awidth]{#3}\\[-3pt] \makebox[\awidth]{#4}}
  }}
\begin{document}

\pagestyle{myheadings}
\date{}
\markboth{\sc Veerman, Peixoto, Rocha, and Sutherland}
         {\sc On Brillouin Zones} 

\title{ON BRILLOUIN ZONES}

\settowidth{\awidth}{\footnotesize Mathematics Departm}
\author{
  \hbox{
    \abox{J. J. P. Veerman}{Physics Department}{UFPE}
	{Recife, Brazil}{veerman@dmat.ufpe.br}
    \abox{M. M. Peixoto}{IMPA}{Rio de Janeiro, Brazil}{\hbox{}}
        {peixoto@impa.br}
    \abox{A. C. Rocha}{Mathematics Dept.}{UFPE}{Recife, Brazil}
        {acr@dmat.ufpe.br}
    \abox{S. Sutherland}{Mathematics Dept.}{SUNY}{Stony Brook, NY, USA}
	{scott@math.sunysb.edu}
 }
}

\maketitle
\SBIMSMark{1998/7}{June 1998}{}
\thispagestyle{empty}


\begin{abstract}
Brillouin zones were introduced by Brillouin~\cite{Br} in the thirties to
describe quantum mechanical properties of crystals, that is, in a lattice in
$\R^n$. They play an important role in solid-state physics.  It was shown by
Bieberbach~\cite{Bi} that Brillouin zones tile the underlying
space and that each zone has the same area.  We generalize
the notion of Brillouin Zones to apply to an arbitrary discrete set in a
proper metric space, and show that analogs of Bieberbach's results
hold in this context.

We then use these ideas to discuss focusing of geodesics in orbifolds of
constant curvature. In the particular case of the Riemann surfaces
$\H^2/\Gamma (k) \,\,(k=2,3, \text{or} 5)$, we explicitly count the number
of geodesics of length $t$ that connect the point $i$ to itself. 
\end{abstract}


\section{Introduction}
\label{app=intro}
\setcounter{figure}{0}
\setcounter{equation}{0}


In solid-state physics, the notion of Brillouin zones is used to describe
the behavior of an electron in a perfect crystal.  In a crystal, the atoms
are often arranged in a lattice; for example, in NaCl, the sodium and
chlorine atoms are arranged along the points of the simple cubic lattice
$\Z^3$.  If we pick a specific atom and call it the origin, its \dfn{first
Brillouin zone} consists of the points in $\R^3$ which are closer to the
origin than to any other element of the lattice.  This same zone can be
constructed as follows: for each element $a$ in the lattice, let $L_{0a}$ be
the perpendicular bisecting plane of the line between $0$ and $a$ (this
plane is called a Bragg plane).  The volume about the origin enclosed by
these intersecting planes is the first Brillouin zone, $b_1(0)$. This
construction also allows us to define the higher Brillouin zones as well: a
point $x$ is in $b_n$ if the line connecting it to the origin crosses
exactly $n-1$ planes $L_{0a}$, counted with multiplicity.

\begin{figure}[htp]
 	\centerline{\psfig{figure=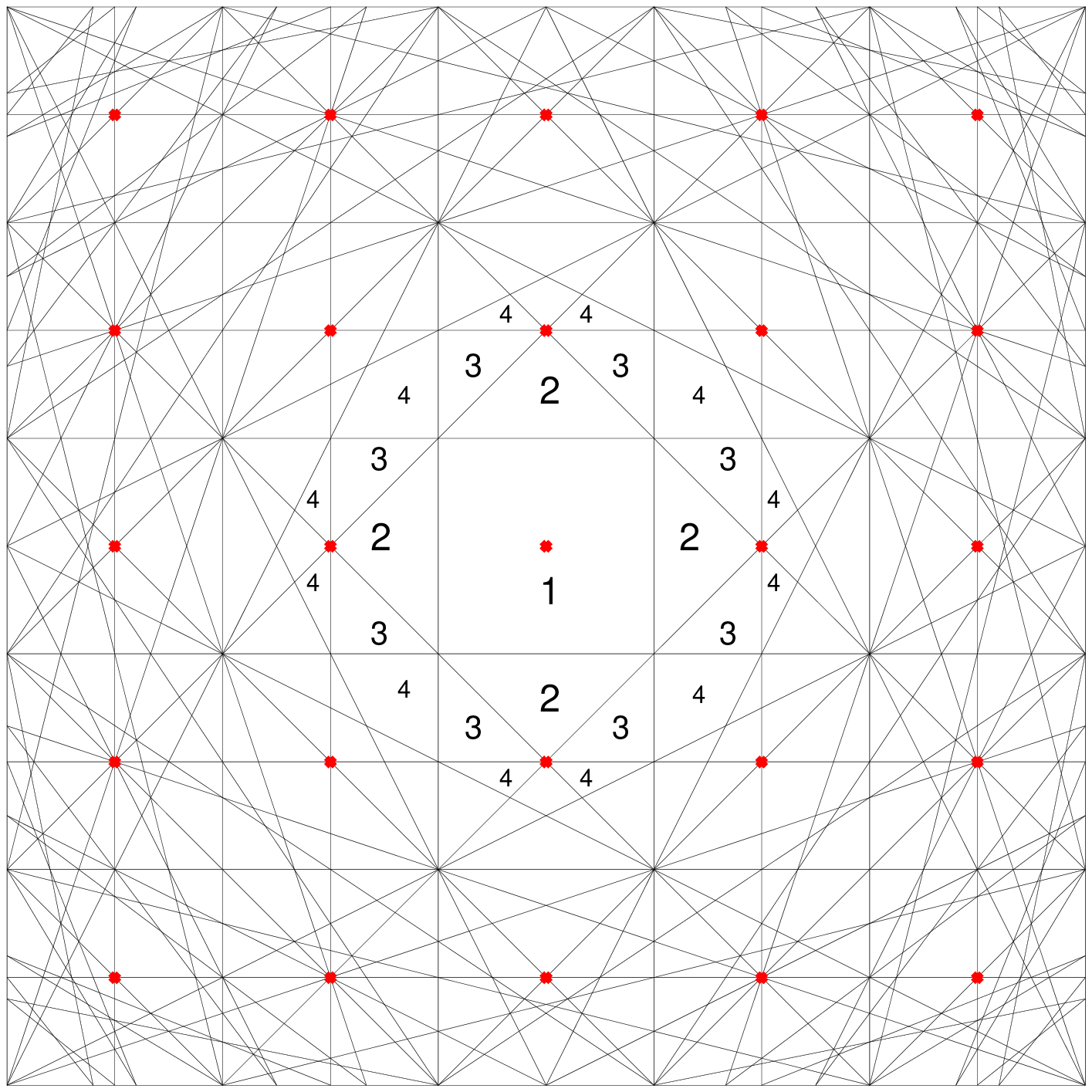,width=.4\hsize} \hfil
                    \psfig{figure=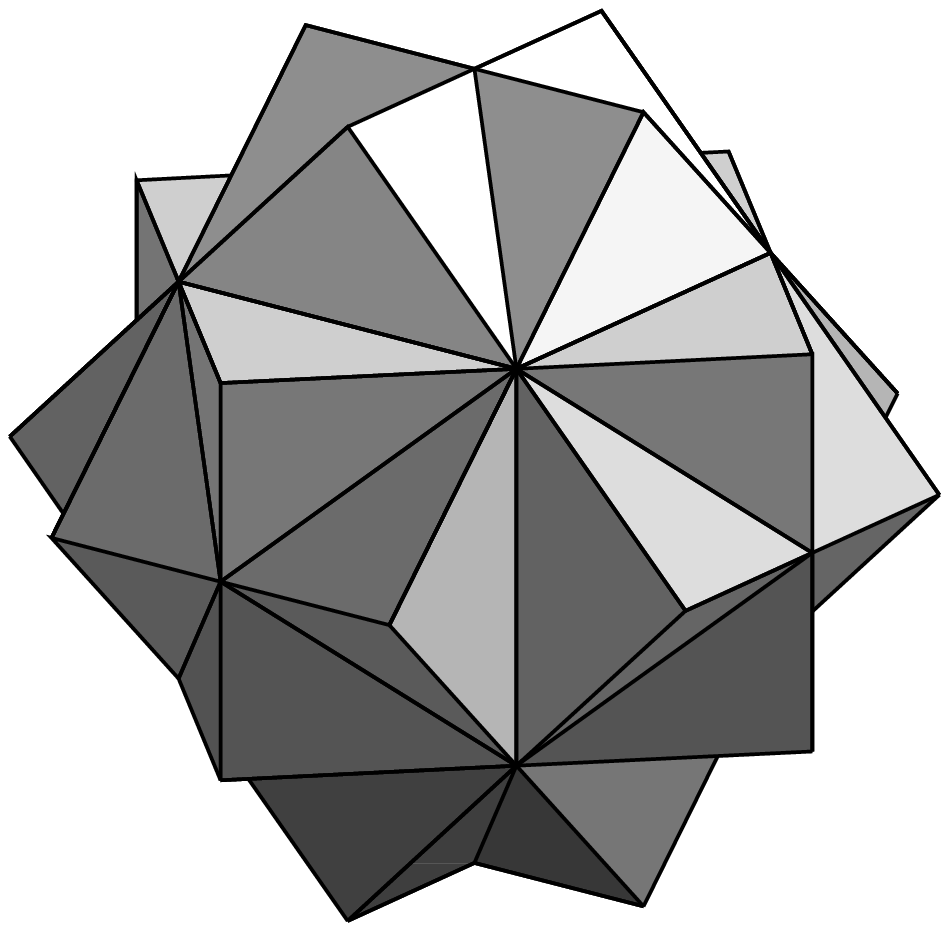,width=.4\hsize} }
 	\Caption{On the right are the Brillouin zones for the lattice $\Z^2$
                    in $\R^2$.  On the left is the outer boundary of the
                    third Brillouin zone for the lattice $\Z^3$ in $\R^3$.}
\end{figure}

This notion was introduced by Brillouin in the 1930s (\cite{Br}), and plays an
important role in solid-state theory (see, for example, \cite{AM} or
\cite{Jo2}, also \cite{Ti}).   The construction which gives rise to Brillouin
zones is not limited to consideration of crystals, however. 
For example, in computational geometry, the notion of the Voronoi cell
corresponds exactly to the first Brillouin zone described above (see
\cite{PS}). We shall also see below how, after suitable generalization, this
construction coincides with the Dirichlet domain of Riemannian geometry, and
in many cases, to the focal decomposition introduced in \cite{Pe1} (see also
\cite{Pe3}).

\bigskip
With some slight hypotheses (see section~\ref{Main Results}), we
generalize the construction of Brillouin zones to any discrete set $S$ in a
path-connected, proper metric space $X$.  We generalize the Bragg planes
above as mediatrices, defined here. 

\begin{defn}
\label{equidistant}
For $a$ and $b$ distinct points in $S$, define the \dfn{mediatrix} (also
called the \dfn{equidistant set} or \dfn{bisector}) $L_{ab}$ of $a$ and $b$
as: 
\bsenn
L_{ab} = \set{x\in X \st d(x,a) = d(x,b)}   .
\esenn
\end{defn}

Now choose a preferred point $x_0$ in $S$, and consider the collection of
mediatrices $\set{L_{x_0,s}}_{s\in S}$.  These partition $X$ into Brillouin
zones as above:  roughly, the $n$-th Brillouin zone $B_n(x_0)$ consists of
those points in $X$ which are accessible from $x_0$ by crossing exactly $n-1$
mediatrices.   (There is some difficulty accounting for multiple crossings---
see definition~\ref{defbrillouin} for a precise statement.)

One basic property of the zones $B_n$ is that they tile the space $X$:
\bsenn
    \Union_{x_i \in S} B_n(x_i)= X \quad\logand\quad
    B_n(x_0)\intersect B_n(x_1) \text{is small}.
\esenn
Here, with some extra hypotheses, ``small'' means of measure zero.
Furthermore, again with some extra hypothesis, each zone $B_n$ has the same
area.  (This property was ``obvious'' to Brillouin).  Both results were
proved by Bieberbach in \cite{Bi} in the case of a lattice in $\R^2$.
Indeed, he proves (as we do) that each zone forms a fundamental set for the
group action of the lattice.  His arguments rely heavily on planar Euclidean
geometry, although he remarks that his considerations work equally well in
$\R^d$ and can be extended to ``groups of motions in non-Euclidean spaces''.
In \cite{G.A.Jones}, Jones proves 
these results for lattices in $\R^d$, as well as giving asymptotics for both
the distance from $B_n$ to the basepoint, and for the number of connected
components of the interior of $B_n$.  In section~\ref{Main Results}, we show
that the tiling result holds for arbitrary discrete sets in a metric space.
If the discrete set is generated by a group of isometries, we show that each
$B_n$ forms a fundamental set, and consequently all have the same area (see
Prop.~\ref{fund_set}).

\bigskip

\medskip
We now discuss the relationship of Brillouin zones and focal decomposition
of Riemannian manifolds.

If $x_1(t)$ and $x_2(t)$ are two solutions of a second order differential
equation with $x_1(0) = x_2(0)$ and there is some $T\ne 0$ so that $x_1(T) =
x_2(T)$, then the trajectories $x_1$ and $x_2$ are said to \dfn{focus} at time
$T$.  One can ask how the number of trajectories which focus varies with the
endpoint $x(T)$--- this gives rise to the
concept of a \dfn{focal decomposition} (originally called a \dfn{sigma
decomposition}).  This concept was introduced in \cite{Pe1} and  
has important applications in physics, for example when computing the
semiclassical quantization using the Feynman path integral method (see
\cite{Pe3}).  There is also a connection with the arithmetic of positive
definite quadratic forms (see \cite{Pe2}, \cite{KP}, and \cite{Pe3}).
Brillouin zones have a similar connection with arithmetic, as can be seen in
section~\ref{examples}.

More specifically, consider the two-point boundary problem
$$\ddot{x} = f(t,x,\dot{x}),  \qquad x(t_0)=x_0, \qquad x(t_1)=x_1,
\qquad x, t, \dot{x}, \ddot{x} \in \R.$$
Associated with this equation, there is a partition of $\R^4$ into sets
$\Sigma_{k}$, where $(x_0,x_1,t_0,t_1)$ is in $\Sigma_{k}$
if there are exactly $k$ solutions which connect $(x_0,t_0)$ to
$(x_1,t_1)$.  This partition is the focal decomposition with respect to the
boundary value problem.  In \cite{PT}, several explicit examples are worked
out, in particular the fundamental example of the pendulum $\ddot{x} = -\sin
x$.  Also, using results of Hironaka (\cite{Hi}) and Hardt (\cite{Ha}), 
the possibility of a general, analytic theory was pointed out.  In
particular, under very general hypotheses, the focal decomposition yields an
analytic Whitney stratification.

Later, in \cite{KP}, the idea of focal decomposition was approached in the
context of geodesics of a Riemannian manifold $M$ (in addition to a
reformulation of the main theorem of \cite{PT}).  Here, one takes a basepoint
$x_0$ in the manifold $M$: two geodesics $\gamma_1$ and $\gamma_2$
\dfn{focus} at some point $y\in M$ if $\gamma_1(T) = y = \gamma_2(T)$.  This
gives rise to a decomposition of the tangent space of $M$ at $x$ into regions
where the same number of geodesics focus.

\medskip
In order to study focusing of geodesics on an orbifold $(M,g)$ with metric
$g$ via Brillouin zones, we do the following.   Choose a base-point $p_0$ in
$M$ and construct  the universal cover $X$, lifting $p_0$ to a point $x_0$ in
$X$. Let $\gamma$ be a smooth curve in $M$ with initial point $p_0$ and
endpoint $p$. Lift $\gamma$ to ${\tilde \gamma}$ in $X$ with initial
point $x_0$. Its endpoint will be some $x \in \pi^{-1}(p)$. The metric
$g$ on $M$ is lifted to a metric ${\tilde g}$ on $X$ by setting
${\tilde g} = \pi^* g$.
Under the above conditions, the group $G$ of deck transformations is 
discontinuous and so $\pi^{-1}(p_0)\subset X$ is a discrete set.
One can ask how many geodesics of length $t$ there are which start at $p_0$
end in $p$, or translated
to $(X,{\tilde \gamma})$, this becomes: {\it How many mediatrices
$L_{x_0,s}$ intersect at $x$, as $s$ ranges over $\pi^{-1}(p_0)$?}

Notice that if the universal cover of $M$ coincides with the tangent space
$TM_x$, the focal decomposition of \cite{KP} and that given by Brillouin zones
will be the same.  If the universal cover and the tangent space are
homeomorphic (as is the case for a manifold of constant negative
curvature), the two decompositions are not identical, but there is a clear
correspondence.  However, if the universal cover of the manifold is not
homeomorphic to the tangent space at the base point, the focal decomposition
and that given by constructing Brillouin zones in the universal cover are
completely different.  For example, let $M$ be $\S^n$, and let $x$ be any
point in it. The focal decomposition with respect to $x$ gives a collection
of nested $n-1$-spheres centered at $x$; on each of these infinitely many
geodesics focus (each sphere is mapped by the exponential to either $x$ or
its antipodal point). Between the spheres are bands in which no focusing
occurs. (See \cite{Pe3}).  However, using the construction outlined in the
previous paragraph gives a very different result.  Since $\S^n$ is simply
connected, it is its own universal cover.  There is only one point in our
discrete set, and so the entire sphere $\S^n$ is in the first zone $B_1$. 

\bigskip
The organization of this paper is as follows.  In section~\ref{Main Results},
we set up the general machinery we need, and prove the main theorems in the
context of a discrete set $S$ in a proper metric space.  

Section~\ref{groups} explores this in the context of manifolds of constant
curvature. The universal cover is $\R^n$, $\S^n$, or $\H^n$, and the  group
$G$ of deck transformations is a discrete group of isometries (see
\cite{Ca}). The discrete set $S$ is the orbit of a point not fixed by any
element of $G$ under this discontinuous group. It is easy to see that the
mediatrices in this case are totally geodesic spaces. From the basic
property explained above, one can deduce that the $n$-th Brillouin zone is a
fundamental domain for the group $G$ in $X$.  

In section~\ref{examples}, we calculate exactly the number of geodesics of
length $t$ that connect the origin to itself in two cases: the flat torus
$\R^2/\Z^2$ and the Riemann surfaces $\H^2/\Gamma(p)$, for $p \in
\set{2,3,5}$.   While these calculations could, of course, be done
independent of our construction, the Brillouin zones help visualize the
process.   

In the final section, we give a nontrivial example in the case of a
non-Riemannian metric, and mention a connection to the question of how many
integer solutions there are to the equation $a^k + b^k = n$, for fixed $k$.

\bigskip
\noindent
{\bf Acknowledgments:} It is a pleasure to acknowledge useful 
conversations with Federico Bo\-netto, Johann Dupont, Irwin Kra, Bernie
Maskit, John Milnor, Chi-Han Sah, and Duncan Sands.  Part of this work was
carried out while Peter Veerman was visiting  the Center for Physics and
Biology at Rockefeller University and the Mathematics Department at SUNY
Stony Brook; the authors are grateful for the hospitality of these
institutions. 


\section{Definitions and main results}
\label{Main Results}
\setcounter{figure}{0}
\setcounter{equation}{0}

\vskip.2in

In this section, we prove that under very general conditions, Brillouin zones
tile (as defined below) the space in which they are defined, 
generalizing an old result of Bieberbach 
\cite{Bi}. With stronger assumptions, we prove that these tiles are
in fact well-behaved sets: they are equal to the closure of their interior. 

\medskip
\noindent{\bf Notation:}
Throughout this paper, we shall assume $X$ is a path connected, proper (see
below) metric space (with metric $d(\cdot,\cdot)$). 
We will make use the following notation: 

\vskip -.5\baselineskip
\begin{itemize} \itemsetup
\item  Write an open $r$-neighborhood of a point $x_0$ as
$\DS{N_r(x_0)= \set{ x\in X \st d(x_0,x)<r} }$. 
\item Define the circumference as
$\DS{ C_r(x_0)= \set{ x\in X \st d(x_0,x)=r} }$. 
\item Their union is the closed disk of radius $r$, denoted by 
$\DS{D_r(x_0)= \set{ x\in X \st d(x_0,x) \leq r}}$.

\end{itemize}

\begin{defn}
A metric space $X$ is \dfn{proper} if the distance function $d(x,\cdot)$
is a proper map for every fixed $x\in X$.  In particular, for every $x\in X$
and $r>0$, the closed ball $D_r(x)$ is compact. Such a metric space is also
sometimes called a \dfn{geometry} (See \cite{Cannon}). 
\end{defn}

\noindent
Note if $X$ is proper, it is locally compact and complete.  The converse
also holds if $X$ is a geodesic metric space (see Thm.~1.10 of
\cite{Gromov}).  The metric spaces considered here need not be geodesic.

\begin{defn}
\label{metr-cons} The space $X$ is called \dfn{metrically consistent} if,
for all $x$ in $X$, $R > r >0$ in $\R$, and for each $a\in C_R(x)$,
there is a $z \in C_r(x)$ satisfying
             $\DS{N_{d(z,a)}(z) \subseteq N_R(x)}$
and $\DS{C_{d(z,a)}(z) \intersect C_R(x) = \set{a}}$.
\end{defn}
 
\noindent
This property is satisfied for any Riemannian metric. 

\smallskip
Note that any mediatrix $L_{a,b}$ separates $X$, that is: $X-L_{ab}$
contains at least two components (one containing the point $a$ and the other
$b$).   

\begin{defn}
\label{min-sepa}
We say that the mediatrix $L_{ab}$ is \dfn{minimally separating} if for any
subset 
$\tilde{L} \subset L_{ab}$ with $\tilde{L} \neq L_{ab}$, the set $X-\tilde{L}$
has one component. 
\end{defn}

\noindent
If a mediatrix $L$ it is minimally separating, then $X-L$ has exactly two
components.   Note that if a separating set $L \subseteq X$ 
contains a non-empty open set $V$, it cannot be minimally
separating. For if $A$ and $B$ are disjoint open sets containing 
$X-L$, then so are $\tilde{A} = A-(A\cap V)$ and
$\tilde{B}=B-(B\cap V)$. Now let $x$ be any point in $V$. Then
it is easy to see that the disjoint open sets $\tilde{A}\cup V$ and
$\tilde{B}$ cover $X-(L-x)$. Thus $L-x$ separates $X$, so $L$ could not have
been minimal.

We define the following sets  
$$L_{0a}^-=\set{ x\in X \st  d(0,x)-d(a,x) < 0 }, \quad
  L_{0a}^+=\set{ x\in X \st  d(0,x)-d(a,x) > 0 }.$$ 

\begin{defn}
\label{topo-trans}
Two minimally separating sets $L_{0a}$ and $L_{0b}$ are \dfn{topologically
transversal} if they are disjoint, or if for each $x\in L_{0a}\cap
L_{0b}$ and every neighborhood $V$ of $x$, the sets 
$L_{0a}^+\cap L_{0b}^+\cap V$,  $L_{0a}^+\cap L_{0b}^-\cap V$, 
$L_{0a}^-\cap L_{0b}^+\cap V$,  $L_{0a}^-\cap L_{0b}^-\cap V$ 
are all nonempty. 
\end{defn}

\smallskip
Usually, there will be a discrete set of points  $S = \set{x_i}_{i\in I}$ in
$X$ which will be of interest. By discrete we mean that any compact subset
of $X$ contains finitely many points of $S$. Note that if 
$\DS{   \liminf_{a,b\in S} d(a,b) > 0 }$, then $S$ is discrete.

\begin{defn}
\label{conditions} 
We say a proper, path connected metric space $X$ is \dfn{Brillouin} if it
satisfies the following conditions:   
\begin{itemize} \itemsetup
\item[1:] $X$ is metrically consistent.
\item[2:] For all $a$, $b$ in $X$, the mediatrices $L_{ab}$ are minimally
separating sets. 
\item[3:] For any three distinct points $0$, $a$, and $b$ in $X$, the
mediatrices $L_{0a}$ and $L_{0b}$ are topologically transversal.

\end{itemize}
\end{defn}

The last two conditions in the above definition may be weakened to apply
only to those mediatrices $L_{ab}$ where $a$ and $b$ in $S$.  In this
case, we will say that \dfn{$X$ is Brillouin over $S$}, if it is not obvious
from the context.

\begin{figure}[ht]
\hbox to \hsize{\hfil
  \hbox{
  \vbox to .6\hsize{\hsize=.4\hsize
	\centerline{\psfig{figure=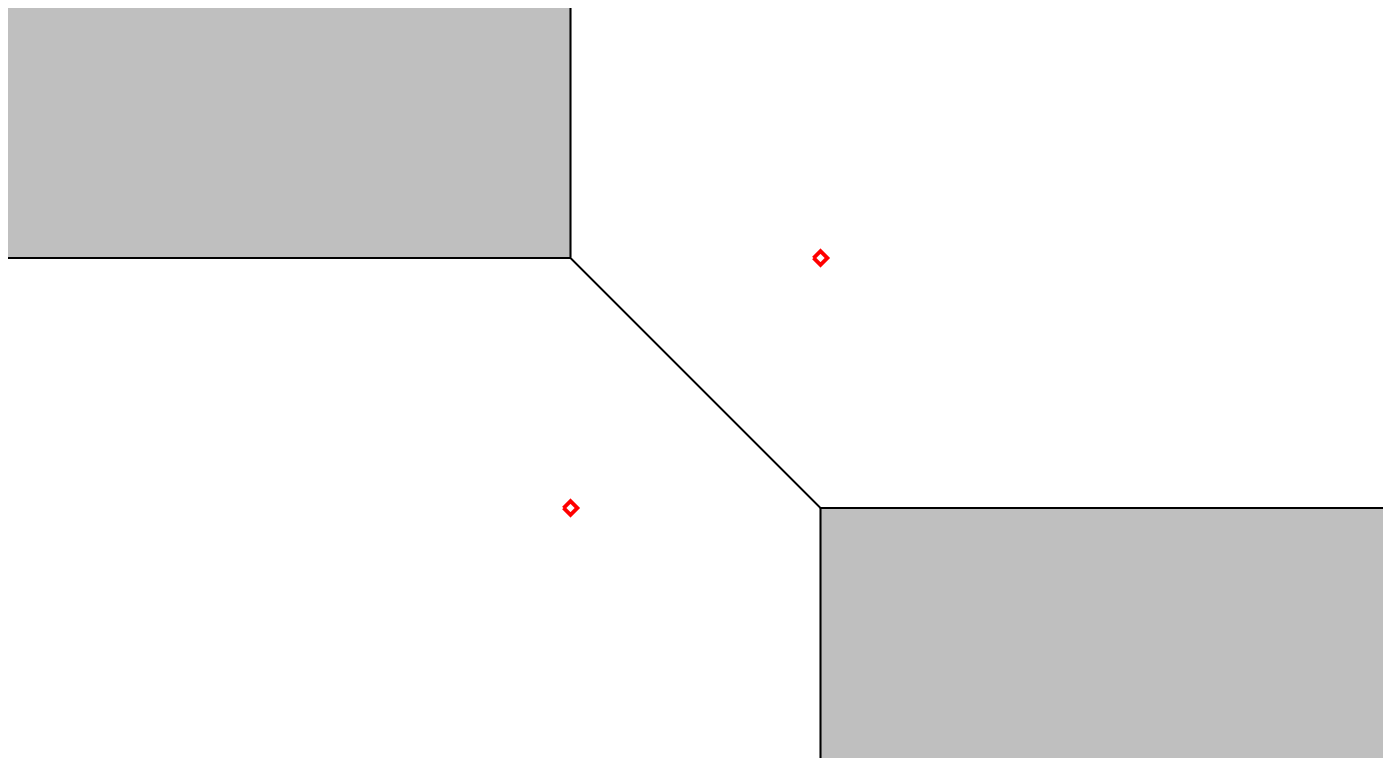,width=.95\hsize}}
	\Caption{\label{fig-man1}The set $L_{(0,0),(a,a)}$ contains
                two quarter-planes.} 
        \vskip 10pt \vfil
	\centerline{\psfig{figure=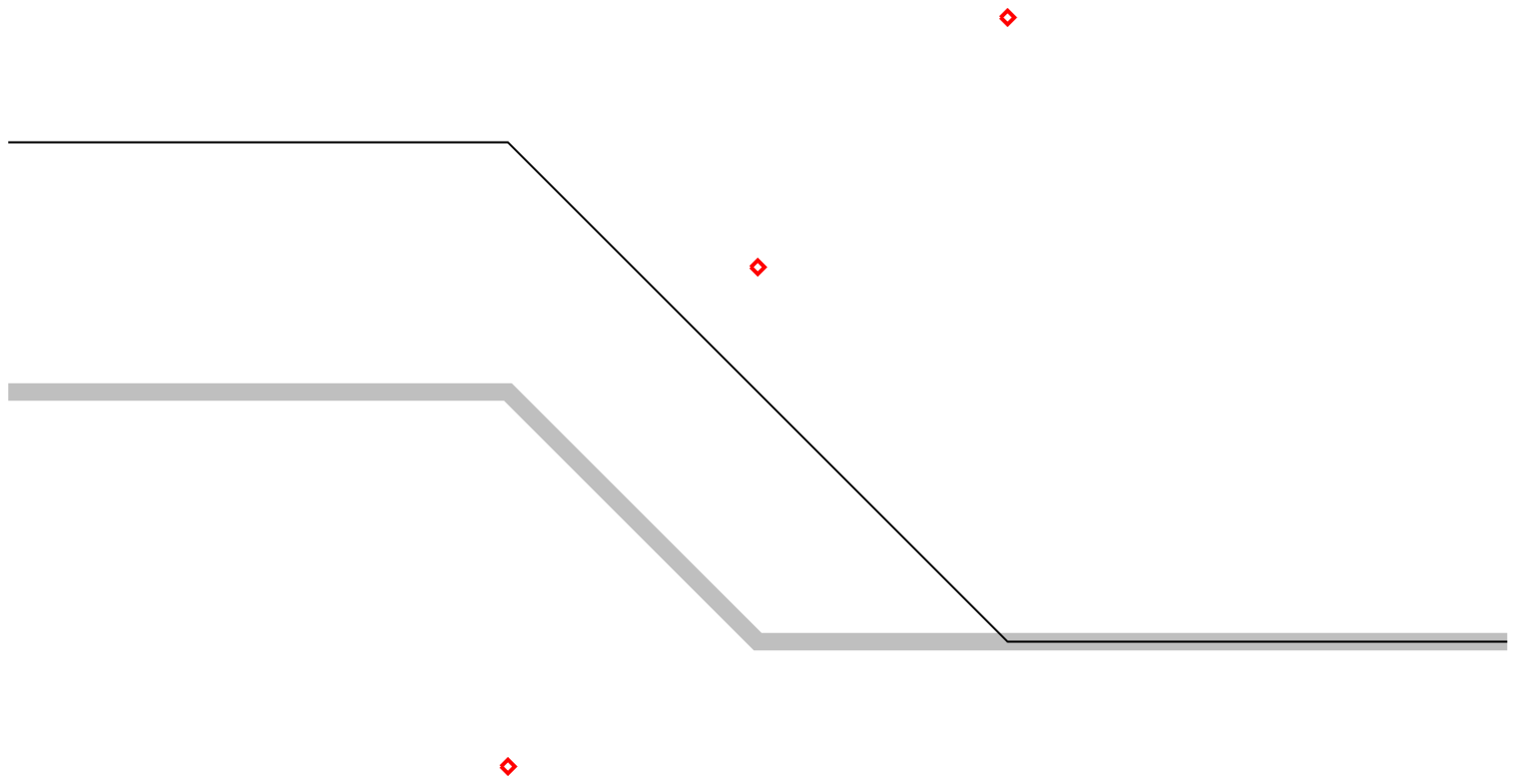,width=.95\hsize}}
	\Caption{\label{fig-man2}
                $L_{(0,0),(4,6)}$ (thin solid line) and $L_{(0,0),(2,4)}$
                (thick grey line) are not transverse.}
  }} 
  \hfil
  \hbox{
  \vbox to .6\hsize{\hsize=.575\hsize \vfil
 	\centerline{\psfig{figure=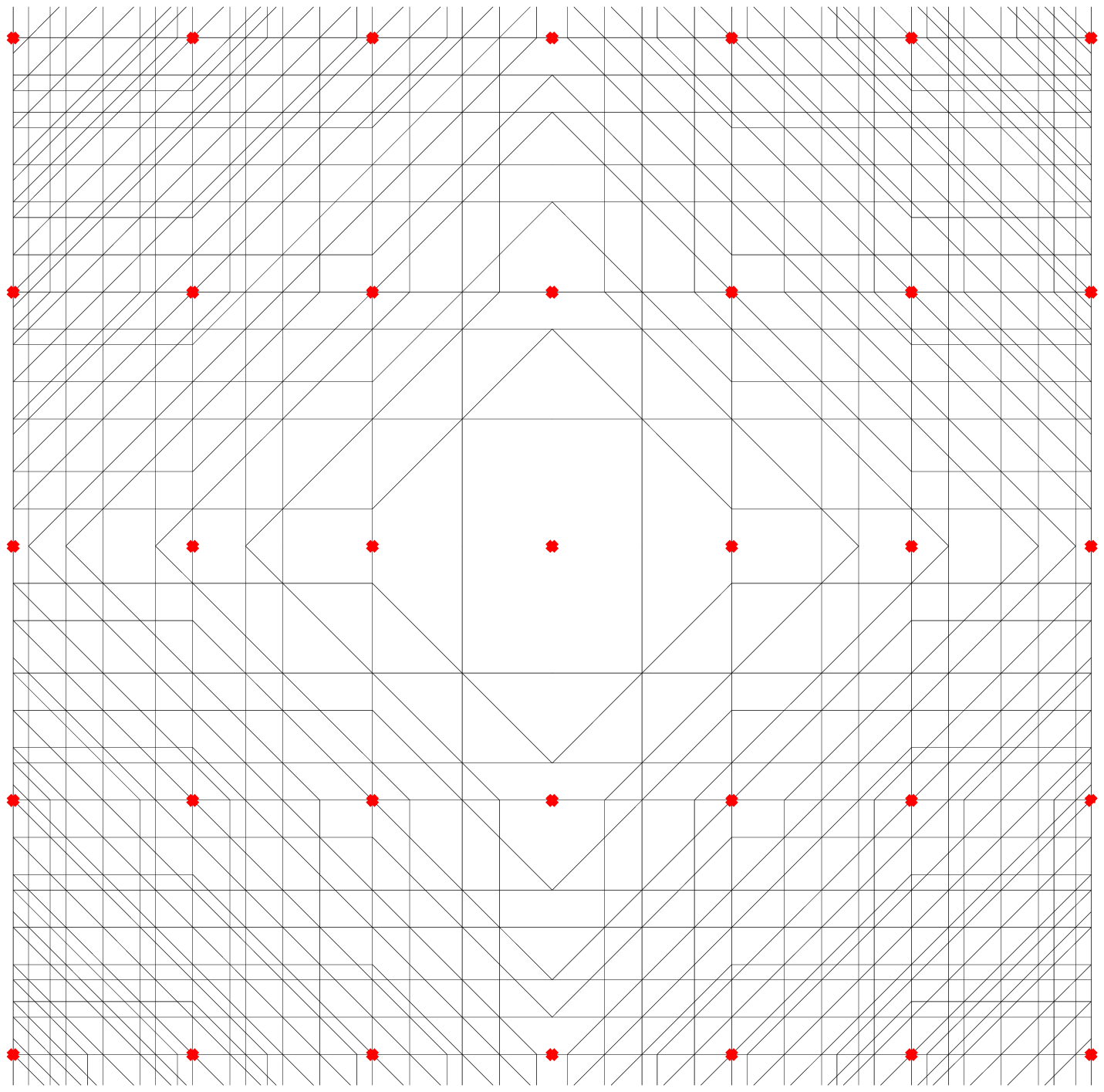,width=.9\hsize}}
	\Caption{\label{fig-man3}
        The mediatrices $L_{0a}$ for $\R^2$ with the Manhattan
	metric and $a$ in the lattice $\set{(m,n\sqrt{2})}$.}
        \vfil
  }} \hfil
}
\end{figure}

\begin{exam}\rm
\label{ex-manhattan}
Equip $\R^2$ with the ``Manhattan metric'', that is,
$d(p,q) = |p_1 - q_1| + |p_2 - q_2|$.  In this metric, a circle $C_r(p)$
is a diamond of side length $r\sqrt{2}$ centered at $p$, so condition
1 is satisfied.  However, condition 2 fails: if the coordinates of a point
$a$ are equal, then $L_{0a}$ consists of a line segment and two
quarter-planes (see Fig.~\ref{fig-man1}).  If the discrete set $S$
contains no such points, we can still run into trouble with topological 
transversality.  For example, the mediatrices $L_{(0,0),(2,4)}$ and
$L_{(0,0),(4,6)}$ both contain the ray $\set{(t,1) \st t \ge 4}$
(Fig.~\ref{fig-man2}). But, if we are careful, we can ensure that the
space is Brillouin over $S$. To achieve this, if $(0,0)$ is the basepoint,
we must have that for all pairs $(a_1,a_2)$ and $(b_1,b_2)$ in $S$, 
$a_1 - a_2 \ne b_1 - b_2$.   For example,  take $S$ to be an irrational
lattice such as $\set{(m,n\sqrt{2}) \st m,n \in \Z}$. It is easy to check
that for all $a,b \in S$, the properties of  definition \ref{conditions} are
true. (From this example, we see that to do well in Manhattan, one should be
carefully irrational.) 
\end{exam}

\medskip
\goodbreak

As mentioned in the introduction, for each $x_0 \in S$, the
mediatrices $L_{x_0 a}$ give a partition of $X$.  Informally, those elements of
the partition which are reached by crossing $n-1$ mediatrices
from $x_0$ form the $n$-th Brillouin zone, $B_n(x_0)$.  This definition is
impractical, in part because a path may cross several mediatrices
simultaneously, or the same mediatrix more than once.  Instead, we will use 
a definition given in terms of the number of elements of $S$ which are
nearest to $x$.  
This definition is equivalent to the informal one 
when $X$ is Brillouin over $S$ (see Prop.~\ref{equiv-defs} below).
We use the notation $\#(S)$ to denote the cardinality of the set $S$.   

\begin{figure}[tbh]
\centerline{\psfig{figure=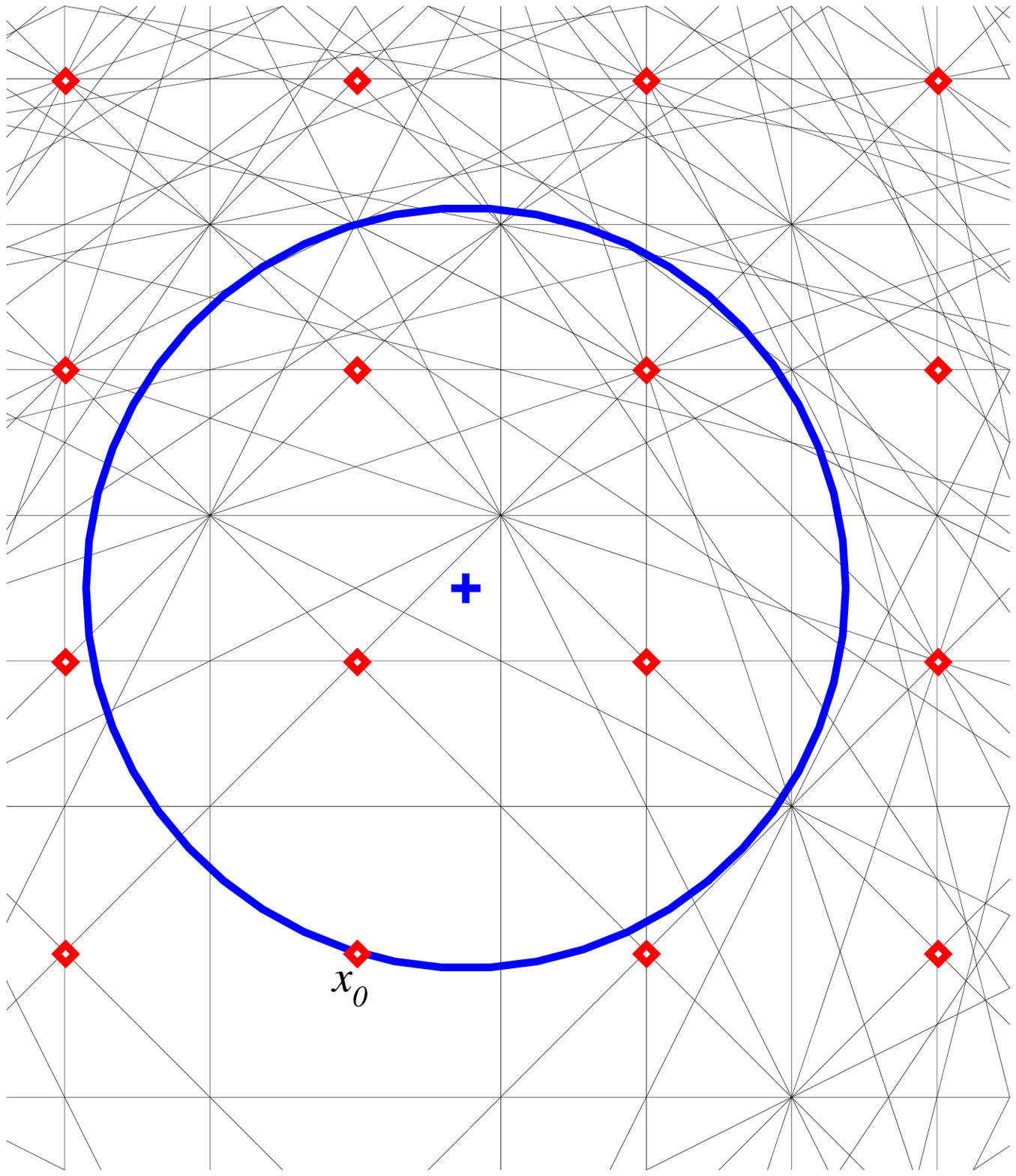,width=.4\hsize} \hfil
            \psfig{figure=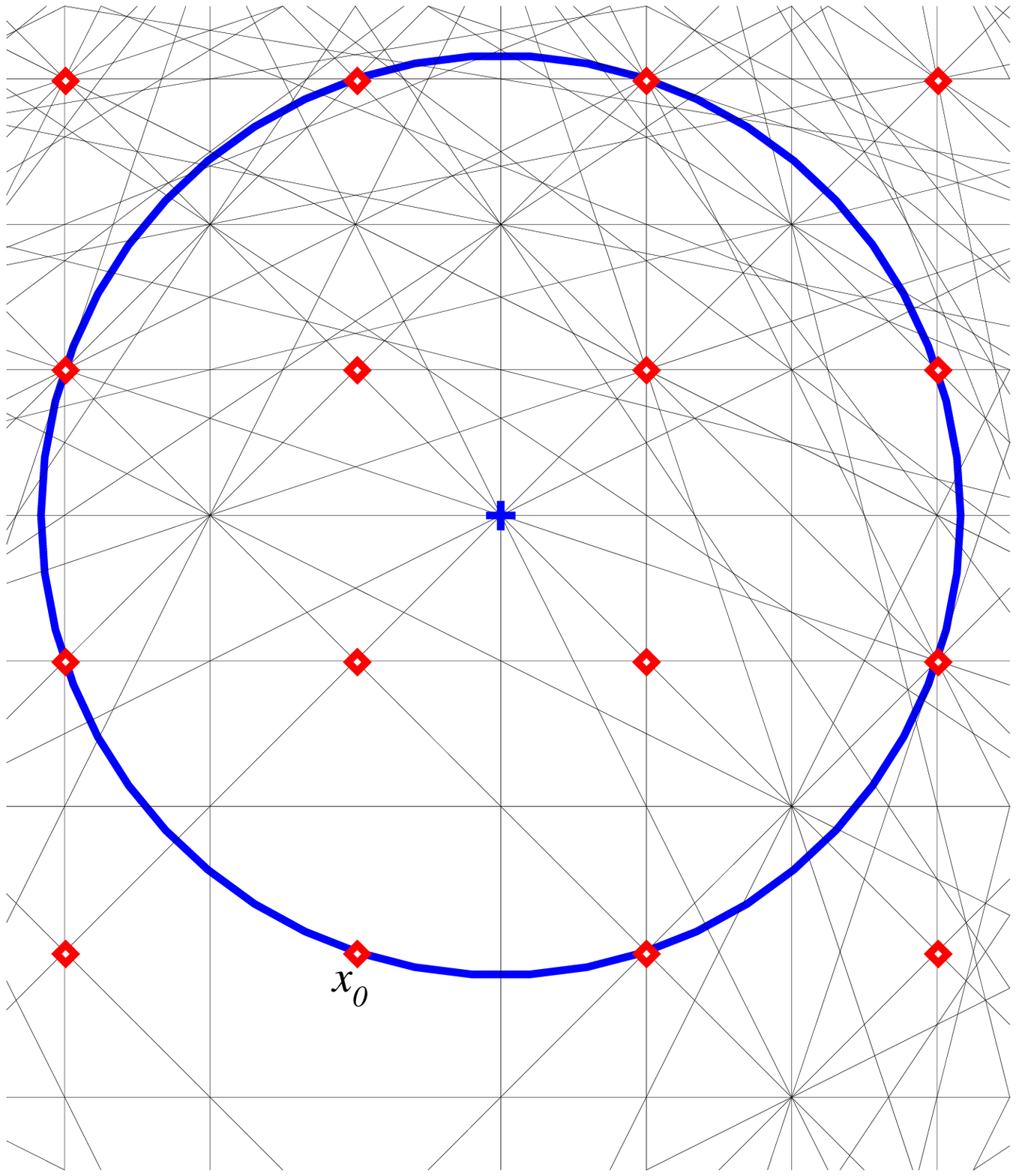,width=.4\hsize} }
\Caption{\label{fig-defbrill} Here we illustrate the definition of the sets
$b_n(x_0)$ and $B_n(x_0)$ for the lattice $\Z^2$ in $\R^2$.  In both
pictures, the circle $\DS{C_{d(x,x_0)}(x)}$ is drawn, and the
basepoint $x_0$ lies in the center of the square at the lower left.
On the left side, the point $x$ (marked by a small cross) lies in
$b_5$, and $\#\left( N_r(x) \intersect S \right) = 4$, while $x_0$ is
the only point of $S$ on the circle.  On the right, we have $m=4$ and
$\ell=8$, so $x$ lies in all of the sets $B_5, B_6, \ldots, B_{12}$.
}
\end{figure}

\begin{defn}
\label{defbrillouin} Let $x\in X$, 
let $n$ be a positive integer, $n \le \#(S)$, and
let $r=d(x,x_0)$. Then define the sets $b_n(x_0)$ and $B_n(x_0)$ as follows.
\begin{itemize} \itemsetup
\item $x \in b_n(x_0)  \iff 
  \#\left( N_r(x) \intersect S \right) = n-1$   \quad and \quad
  $C_r(x) \intersect S = \set{x_0}$.
\item $x \in B_n(x_0) \iff
  \#\left( N_r(x) \intersect S \right) = m$   \quad and \quad
  $\#\left( C_r(x) \intersect S \right) = \ell$,  \quad \mbox{where
  $l,m \in \Z^+$} with \mbox{$m +1 \leq n \leq m +\ell$}.
\end{itemize}
\end{defn}
  
Here the point $x_0$ is called the \dfn{base point}, and the set $B_n(x_0)$
is the 
$n$-th Brillouin zone with base point $x_0$. Note that in the second part, 
if $m=n-1$ and
$\ell=1$, then $x \in b_n(x_0)$. So $b_n(x_0) \subseteq B_n(x_0)$. Note
also that the complement of $b_n(x_0)$ in $B_n(x_0)$ consists of subsets
of mediatrices (see Def.~\ref{equidistant}). 
Note also that $b_n(x_0)$ is open and that
$B_n(x_0)$ is closed. Finally, observe that for fixed $x_0$ the sets
$b_n(x_0)$ are disjoint, but the sets $B_n(x_0)$ are not.  
In what follows it will be proved that $\set{B_n(x_0)}_{n>0}$
cover the space $X$. Thus we can assign to
each point $x$ its Brillouin index as the the largest $n$ for 
which $x\in B_n(x_0)$. This definition was first given in
\cite{Pe1}.

The following lemma, which follows immediately from Def.~\ref{defbrillouin},
explains a basic feature of the zones, namely that they are concentric in in
a weak sense.  This property is also apparent from the figures.

\noindent
\begin{lem}\label{Bn-concentric}
Any continuous path from $x_0$ to $B_n(x_0)$ intersects $B_{n-1}(x_0)$.
\end{lem}


\medskip
The Brillouin zones actually form a covering of $X$ by
non-overlapping closed sets in various ways. This is proved
in parts. The next two results assert that the zones $B$ cover
$X$, but the zones $b$ do not. The first of these  is an immediate
consequence of the definitions. The second is more surprising and leads
to corollary \ref{fundamental}. The fact that the $B_i(x_n)$ 
is the closure of $b_i(x_n)$ (and thus that the interiors do 
not overlap) is proved in proposition \ref{lemma3}. We
will use the word ``tiling'' for a covering by non-overlapping
closed sets.

\noindent
\begin{lem} 
\label{tiling} For fixed n the Brillouin zones tile $X$ in the
following sense: 
\bsenn
	\Union_i B_i(x_n) = X   \quad\logand\quad 
  b_i(x_n) \intersect b_j(x_n) = \emptyset \quad \text{if} i\neq j.
\esenn
\end{lem}
\begin{figure}[tbh]
\centerline{\psfig{figure=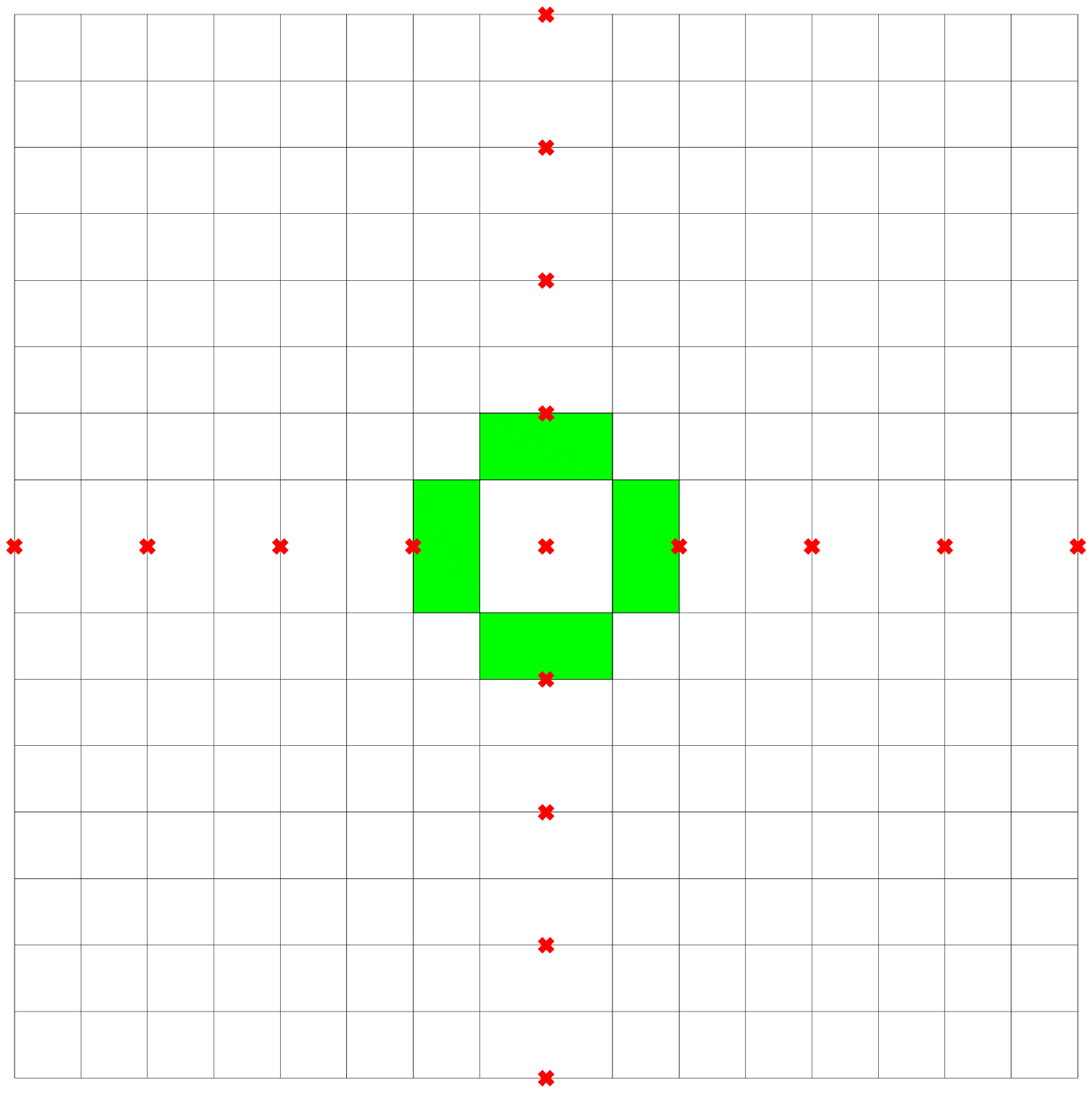,width=.25\hsize} \hfil
            \psfig{figure=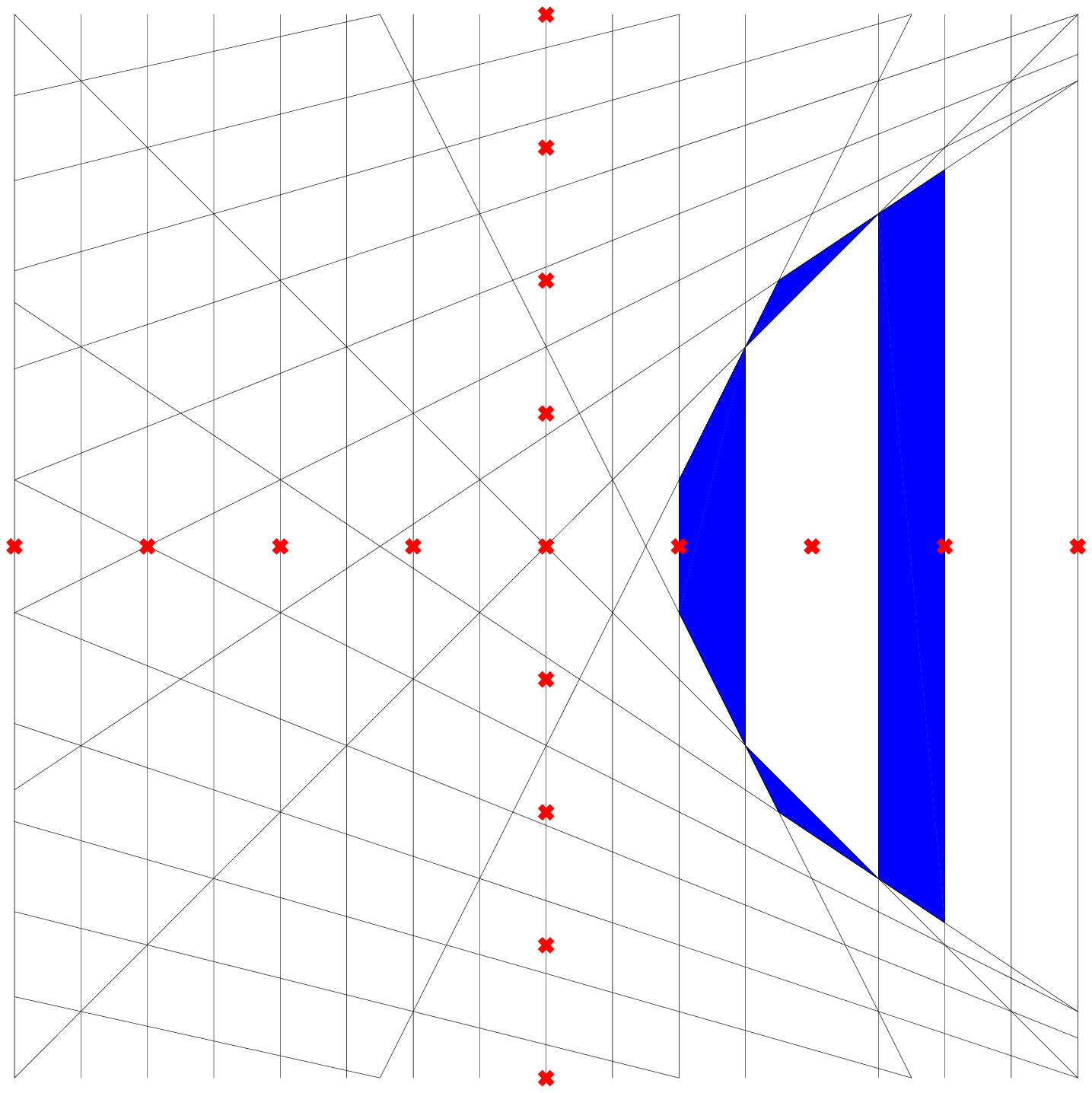,width=.25\hsize} \hfil
            \psfig{figure=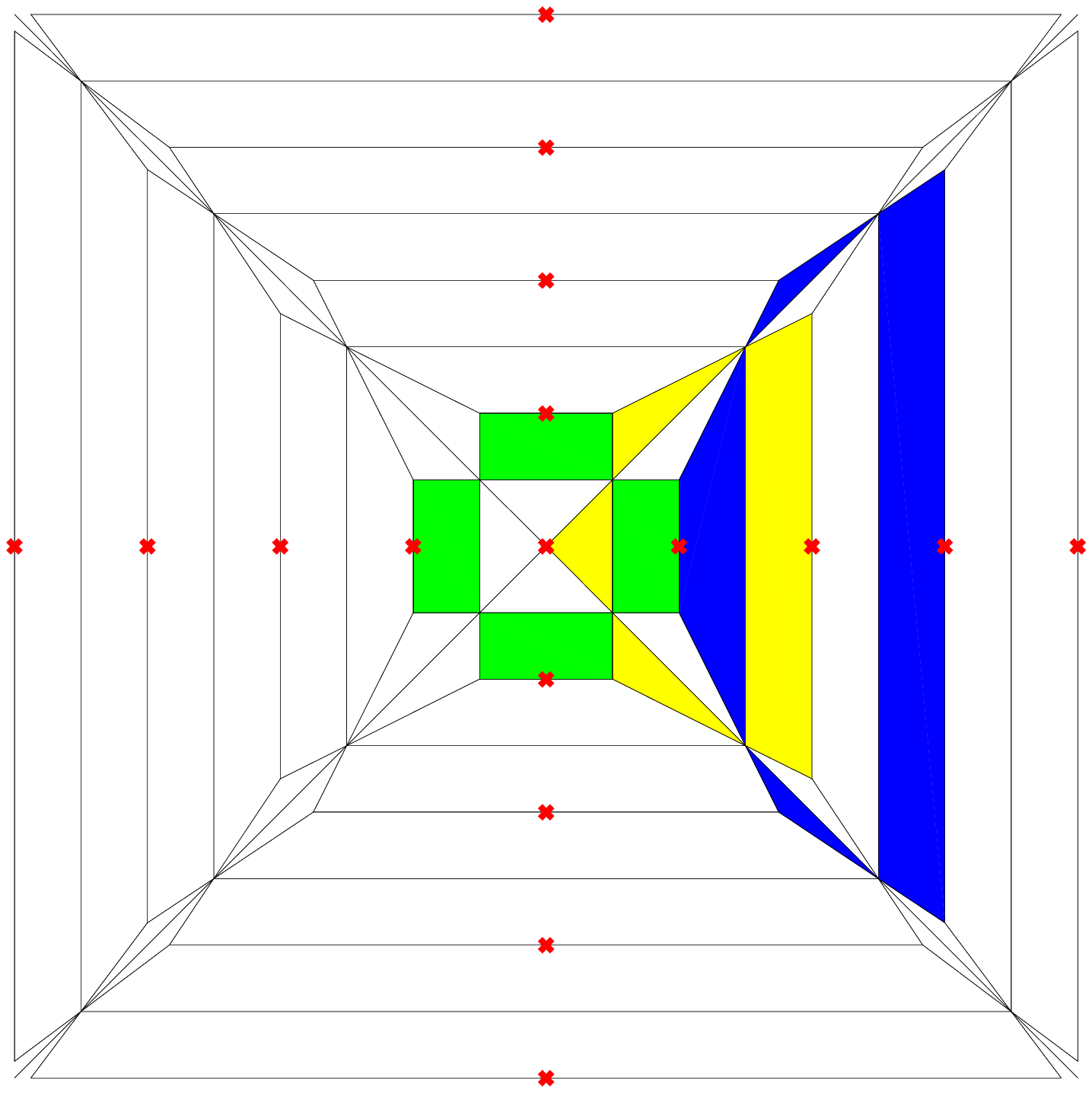,width=.25\hsize} }
\Caption{\label{fig-tiling} This example illustrates Lemma~\ref{tiling} and
Thm.~\ref{bigtheorem}.  Let $S$ be the discrete set $\set{(m,0)}\union
\set{(0,n)}, m,n\in\Z$ in the Euclidean plane. On the left is the
tiling given by $B_i(0,0)$ and in the middle is the tiling by
$B_i(2,0)$. In both cases, $b_2$ is shaded.  On the right is the
tiling given by $B_2(x_i)$ as in Thm.~\ref{bigtheorem}.  The sets
$b_2(0,0)$, $b_2(1,0)$, and $b_2(2,0)$ have been shaded. Note that this
$S$ does not correspond to a group, nor does it satisfy the hypotheses
of Prop.~\ref{fund_set}, because there are no  isometries
which permute $S$ and do not fix the origin.}
\end{figure}

\noindent
\begin{theo}
\label{bigtheorem}
Let $X$ be a proper, path connected metric space let $S = \set{x_i}_{i\in I}$
be a discrete set.  Then, for fixed $n \le \#(S)$, the sets
$\set{B_{n_0}(x_i)}_{i\in I}$ tile $X$ in the following sense:
\bsenn
  \Union_i B_{n}(x_i) = X  \quad\logand\quad
  b_{n}(x_i) \intersect b_{n}(x_j) = \emptyset \quad\text{if} i\neq j.
\esenn
\end{theo}

\proof 
First, we show that for any fixed $n>0$ and each $x \in X$, there is an $x_i
\in S$ with $x \in B_{n}(x_i)$.   Re-index $S$ so that if 
$S =\set{x_1, x_2, x_3,\ldots }$ and  
$i < j$, then $d(x,x_i) \le d(x,x_j)$.
This can be done; since $S$ is a discrete subset and closed
balls $D_c(x_i)$ are compact, 
the subsets of $S$ with $d(x,x_i) \leq c$ are all
finite.  Let $r_i = d(x,x_i)$. We will show that $x \in B_{n}(x_{n})$.

Note that $r_{n} \geq r_{n-1}$. Suppose first that $r_{n} > r_{n-1}$, then 
$N_{r_{n}}(x) \intersect S$
contains exactly ${n}-1$ points, and $x_{n} \in C_{r_{n}}(x) \intersect S$.
Thus $x \in B_{n}(x_{n})$. Note that if $r_{n+1} > r_n$, then we
would have $x \in b_{n}(x_{n}) \subset B_{n}(x_{n})$.

If, on the other hand, $r_{n} = r_{n-1}$,
then there is a $k>0$ so that $r_{n} = r_{n-1} = \ldots = r_{n-k}$, and so
$\#\left( N_{r_{n}}(x) \intersect S \right) = {n} -k-1 \leq {n}-1$.  But then
$\#\left(C_{r_{n}}(x)\intersect S\right) \ge k+1$, and hence 
$x\in B_{n}(x_{n})$ as desired. 

For the second part, we show that $b_{n}(x_i)\intersect b_{n}(x_j) =
\emptyset$.  If not, then there is a point $x$ in their intersection.
If $r_i = r_j$, then $x_i = x_j$, because by the definition of
$b_n(x_k)$, $\set{x_k} = C_{r_k}(x) \intersect S$.  If not, then 
$r_i < r_j$.  In this case, 
$x_i \in D_{r_i}(x) \subset N_{r_j}(x)$ .  Thus, since 
$\#\left( N_{r_i}(x)\intersect S \right) = n-1$, $N_{r_j}(x)$ must
contain at least $n$ points of $S$, a contradiction.
\QED

The next result indicates how this notion of tiling is related to the
notion of fundamental domain. 

\noindent
\begin{prop}
\label{fund_set}
Let $S$ be a discrete set in a metric space $X$ as in theorem 
\ref{bigtheorem}. Suppose that for each $x_i$ in $S$ there is an
isometry $g_i$ of $X$  such that $g_i(x_0) = x_i$, $g_i$ permutes $S$
and the only $g_i$ which leaves $x_0$ fixed is the identity.
Then there is a set $F$ (the fundamental set), satisfying:
\benn
b_n(x_0) \subseteq & F & \subseteq B_n(x_0) \quad \with \\[.5\baselineskip]
\Union_i g_i(F) = X  & \logand &  g_i(F) \intersect g_j(F) = \emptyset 
\; (i\neq j)   .
\eenn
\end{prop}

\proof Suppose that $x \in b_n(x_0)$. From definition \ref{defbrillouin}
and the fact that the $g_i$ are isometries, we see that this is equivalent
to $g_i(x) \in b_n(x_i)$. Thus $g_i(b_n(x_0)) = b_n(x_i)$. Now apply
theorem \ref{bigtheorem}. A similar reasoning proves the statement
for $B_n(x_0)$.
\QED

\remark The fundamental set $F$ is not necessarily connected.
Also, note that it follows from this proposition that $B_i(x_0)$ is
scissors congruent to $B_j(x_0)$ (see \cite{Sah} for a discussion of
scissors congruence).  In particular, this implies immediately that the
$B_i$ all have the same area.  Note that this result does not hold if $S$ is
not generated by a group of isometries.  For example, the ``face-centered
cubic'' structures common to many crystals do not satisfy this, and the
corresponding $B_n$ have different volumes, although they still tile the
space.  See also figure~\ref{fig-tiling}.

\bigskip
The preceding results are valid without the transversality
conditions on the metric, as the reader can check. 
We will now use the additional properties of the metric to derive
a regularity result.

\begin{prop}
\label{lemma3} If $X$ is Brillouin (over S), then 
$B_n(x_i)$ is the closure of $b_n(x_i)$.
\end{prop}

\proof
Let $x\in B_n(x_0)$. If $x\in b_n(x_0)$, we are done.
So suppose $x\not\in b_n(x_0)$, and then by Def.~\ref{defbrillouin} (with 
$r=d(x,x_0)$), $N_r(x)$ contains $m \leq n-1$ points of $S$ and $C_r(x)$ 
contains $\ell \geq n-m \with \ell \geq 2$ points of $S$, including $x_0$. 
Thus $x$ lies in the intersection of $\ell -1 \geq 1$ mediatrices 
$\set{L_{x_0,x_i}}_{i=1}^{\ell -1}.$ 

Suppose $\ell = 2$, so $x$ lies on a single mediatrix $L_1$. Then by
condition~2, $x$ is in the closure of $L_{1}^-$ 
and $L_{1}^+$. Thus if $x' \in L_{1}^-$ and close enough to $x$, then
$x' \in b_{m+1}(x_0)$. If $x' \in L_{1}^+$ and sufficiently close to
$x$, then $x'\in b_{m+2}(x_0)$. From our choices, we see that
$m=n-1$ or $m=n-2$.  In either case, we are done.

The proof for $\ell >2$ is slightly more involved (we have to use
transversality), but the idea is the same. The point $x$ lies in the
intersection of at least two mediatrices.  Essentially, we need show
that the mediatrices disentangle in every neighborhood of $x$.
We will prove that for all $j \in \set{1,\ldots, \ell}$, 
$x \in \closure{b_{m+j}(x_0)}$. Note that $n-\ell \leq m \leq n-1$,
and so showing this is sufficient.

Again, since mediatrices are minimally separating, we have that for
all $i\in \set{1,\ldots ,\ell -1}$ there exist open sets $L_{i}^-$ and 
$L_{i}^+$ in any neighborhood of $x$ for which
	$$d(x,x_0) < d(x,x_i) \text{for} x\in L_{i}^-
	   \quad \text{and} \quad
	  d(x,x_0) > d(x,x_i) \text{for} x\in L_{i}^+.
	$$
Noting that $x_i\in C_r(x),$ then from condition~1 of
Def.~\ref{defbrillouin}, we have that for any small $\rho>0$,  there is a point
$z^- \in C_{\rho}(x)$ for which $d(z^-,x_0) < d(z^-,x_i)$ for all $i\in
\set{1,\ldots, \ell}$. Thus $\bigcap_i L_{i}^-$ is nonempty,
and is contained in $b_{m+1}(x_0)$.
By the same principle, we can choose a point $z^+ \in C_r(x)$ for
which $d(z^+,x_0) > d(z^+,x_i)$.  As a consequence, $\bigcap_i
L_{i}^+$ is nonempty, and is a subset of $b_{m+\ell}(x_0)$.  

Now, $b_{m+1}(x_0)$ and $b_{m+\ell}(x_0)$ are disjoint. 
If $\ell>1$, their complement in a neighborhood of $x$ must contain at
least two mediatrices. By transversality (condition 3), these cannot
coincide. Thus $V$ must be intersected by some $b_{m+i}(x_0)$ for 
$i\in \set{1,\ldots \ell}$.  By induction one finishes the argument.
\QED

\remark The references to mediatrices in the above proof may make it
appear that we are relying on the informal definition of $B_n$ rather than 
Definition~\ref{defbrillouin}.  In fact, we are only using the fact that
$B_n \setminus b_n$ consists of subsets of mediatrices, which was noted
following Def.~\ref{defbrillouin}.  
In fact, for Brillouin spaces, the informal definition of $b_n$ by
mediatrices and Definition~\ref{defbrillouin} are equivalent, as the
following proposition shows.


\noindent
\begin{prop}\label{equiv-defs}
Let $X$ be Brillouin over $S$.  Then there is a path $\gamma$ from
$x$ to $x_0$ which crosses exactly $n-1$ mediatrices with no multiple
crossings if and only if $x \in b_n(x_0)$.  By ``no multiple crossings'' we
mean that if $L_i$ and $L_j$ are distinct mediatrices, $\gamma$ intersects
$L_i$ in at most a single point, and 
$\gamma \intersect L_i \ne \gamma \intersect L_j$. 
\end{prop}

\proof
First, suppose there is such a path $\gamma$ from $x_0$ to $x$.  As
long as $\gamma(t)$ crosses no mediatrices, the number of elements of
$S$ contained in $D_r(\gamma(t))$ remains constant.  If this number is
$i$, then $\gamma(t) \in b_i(x_0)$.  Each time $\gamma(t)$ crosses a single
mediatrix, this number increases by $1$.  Thus, after crossing $n-1$
mediatrices, $\gamma(t)$ will be in $b_n(x_0)$.

Now, let $x\in b_n(x_0)$.  We must construct a path from $x$ to $x_0$
which crosses exactly $n-1$ mediatrices.  Let $b_n^*$ denote the
connected component of $b_n(x_0)$ which contains $x$, and choose a
point $y$ in $B_{n-1}(x_0)  \intersect \partial{b_n^*}$.  Such a point
always exists by Lemma~\ref{Bn-concentric} and Prop.~\ref{lemma3}.  In
fact, because $\partial{b_n}$ consists of subsets of mediatrices,
topological transversality of mediatrices allows us to choose $y$ so that it
lies on a single mediatrix.

There certainly is a path connecting $x$ to $y$ which crosses no
mediatrices; by similar reasoning, we can choose $x' \in b_{n-1}(x_0)$ and a
path from $x'$ to $y$ which crosses no mediatrices.  We have now constructed
a path from $x\in b_n$ to $x' \in b_{n-1}$ which crosses exactly one
mediatrix--- repeating the argument $n-2$ times finishes the proof.  
\QED

\section{Brillouin zones in spaces of constant curvature}
\label{groups}
\setcounter{figure}{0}
\setcounter{equation}{0}

\bigskip

In this section $X$ will be assumed to be 
one of $\R^n$, $\S^n$, or $\H^n$, all equipped with the standard
metric,  and let $G$ be a discontinuous group of
isometries of $X$. Denote the quotient $X/G$ with the induced metric
by $(M,g)$. Then the construction of lifting to the universal cover, as 
outlined in the introduction, applies naturally to $(M,g)$. In this section we
describe focusing of geodesics  in $(M,g)$ by Brillouin zones in $X$. The
discrete set $S$ is given by the orbit of a chosen point in $X$ (which we will
call the origin) under the group of deck-transformations $G$. The fact that
the Brillouin zones are fundamental domains is now a direct corollary of
proposition \ref{fund_set}.  

\medskip
The regularity conditions of Def.~\ref{conditions}
are easily verified in the present context. We do this first.

\begin{lem}
\label{tot-geod}  If $X$ is either $\R^n$, $\S^n$, or $\H^n$, then a mediatrix
$L_{ab}$ in $X$ is an $(n-1)$-dimensional, totally geodesic subspace
consisting of one component, and $X-L_{ab}$ has two components. 
\end{lem}

\proof
This is easy to see if we change coordinates by an isometry of $X$, putting
$a$ and $b$ in a convenient position, say as $x$ and $-x$.  The mediatrix
$L_{x,-x}$ is easily seen to satisfy the conditions (in the case of $\S^n$,
it is the equator, and for the others, it is a hyperplane). The conclusion
follows. 
\QED

\begin{prop}
\label{regularity}
All such spaces $X$ are Brillouin (see definition \ref{conditions}). 
\end{prop}

\proof As remarked before, the first condition is satisfied for any Riemannian 
metric. The second condition is also easy. It suffices to observe
that
the subspaces of Lemma~\ref{tot-geod} are minimally separating.

To prove the topological transversality condition, note that lemma
\ref{tot-geod} 
implies that if $L_{0a}$ and $L_{0b}$ coincide in some open set, then they
are equal. One easily sees that then we must have $L_{0a}=L_{0b}=L_{ab}=L$.
Applying the second part of the lemma to $L_{0a}$, $L_{0b}$, and $L_{ab}$,
we see 
that they must separate $X$ in at least 3 components, namely one for 
each point $a$, $b$, and $0$. This leads to a contradiction.
\QED

\remark Note that in fact, the mediatrices are transversal in the
usual sense. If $L_{0a}$ and $L_{0b}$ coincide in an open set, then their
tangent spaces also coincide at some point. Uniqueness of solutions of 
second order differential equations then implies $L_{0a} = L_{0b}$.

\bigskip
Recall that a metric space $X$ is called {\it rigid\/} if
the only isometry which fixes each point of a
nonempty open subset of $X$ is the identity.
It is not hard to see that $\S^n$, $\H^n$, and $\R^n$ are rigid spaces. 
See \cite{Ra} for more details of rigid metric spaces and for
the proof of the following result. Recall that 
the stabilizer in $G$ of a point $x\in X$ consists of those elements of
$G$ that fix $x$.

\begin{prop}
\label{trivial-stab}
Let $G$ be a discontinuous group of isometries
of a rigid metric space $X$. Then there exists a point $y$ of $X$
whose stabilizer $G_y$ is consists of the identity.
\end{prop}
 
We now return to Brillouin zones as defined in the last section. Recall
that $G$ is a group of isometries of $X$ that acts discontinuously
on points in $X$. Let $x_0$ be a point in $X$ whose stabilizer under
the action of $G$ is trivial. Let $x\in X$ be a point in the orbit
of $x_0$ under $G$, that is to say, $x=\gamma (x_0)$ for some 
$\gamma \in G -\set{e}$, where $e$ is the identity of $G$. Let 
$[x_0,x]$ be a geodesic segment of minimal length whose endpoints are $x_0$
and $x$. Then $B_n(x_0)$, the $n$-th Brillouin zone relative to $x_0$, is
the set of points in $X$ such that the segment $[x_0,x]$ intercepts 
exactly $n-1$ mediatrices $L_{x_0,y}$, where $y$ is in the orbit of $x_0$
under the group $G$. Proposition~\ref{fund_set} immediately implies the most
important fact about Brillouin zones in this setting. 

\begin{figure}[ht]
  \centerline{\psfig{figure=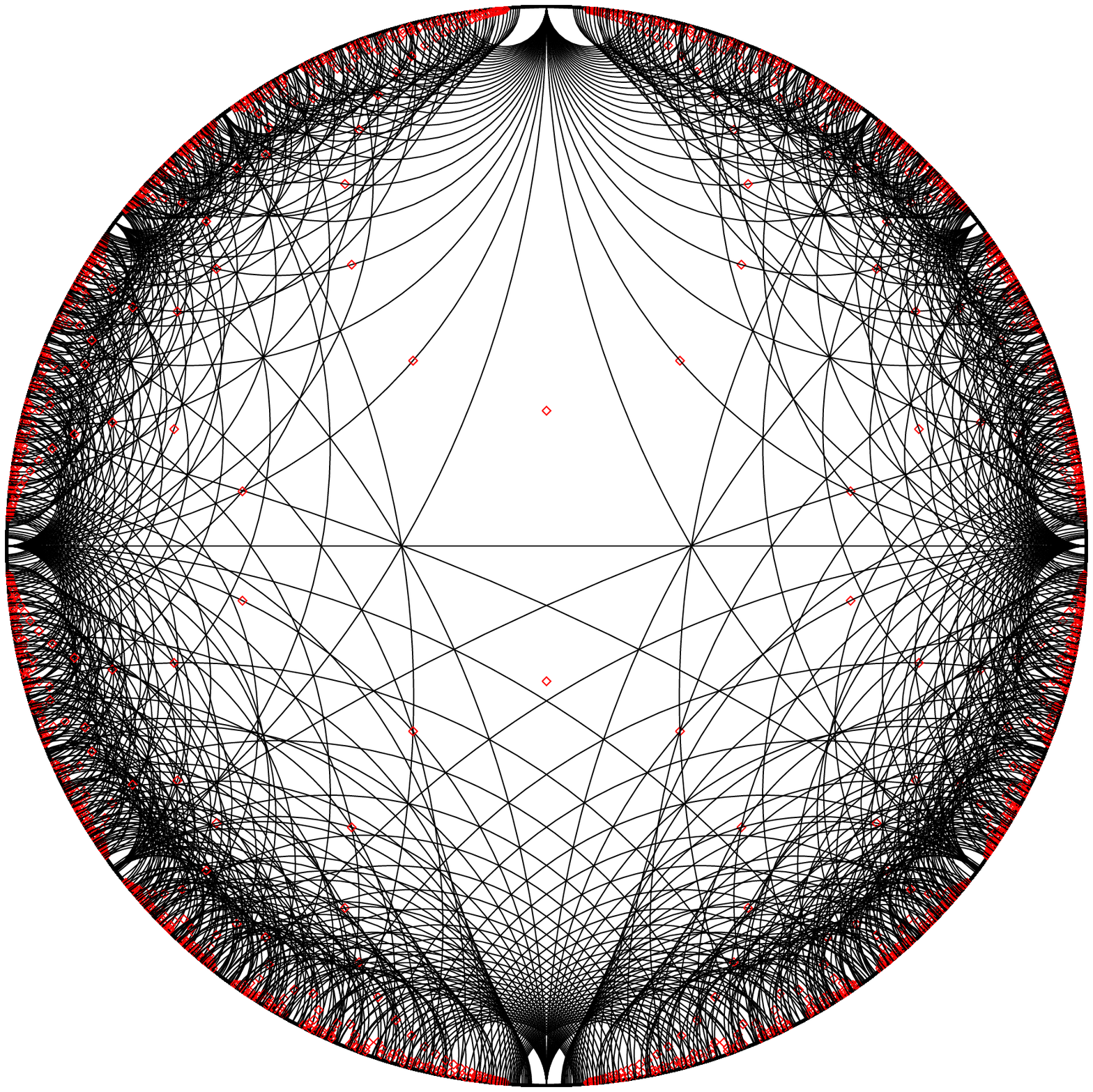,width=.4\hsize}
	     \qquad
             \psfig{figure=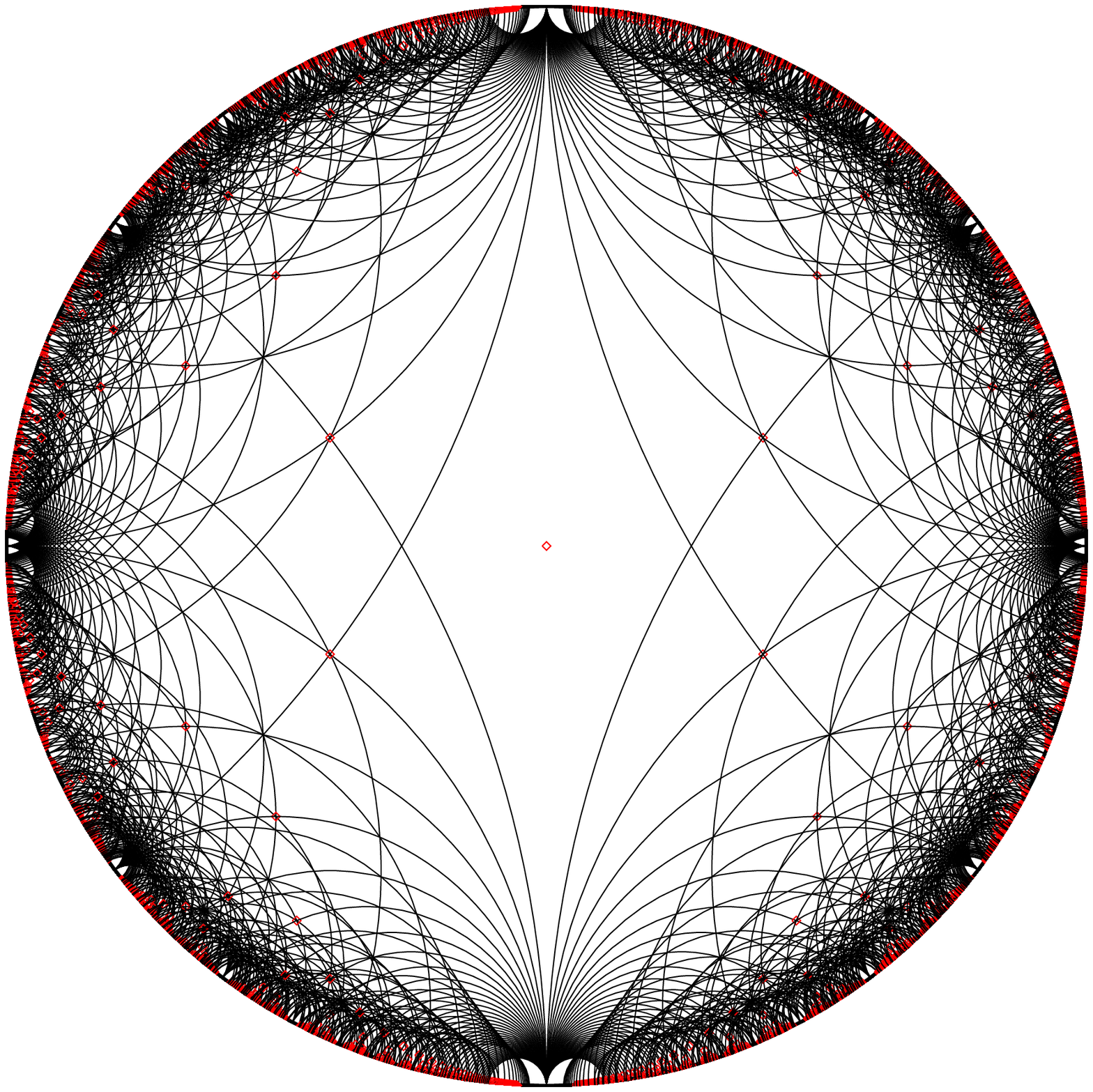,width=.4\hsize}}
 \Caption{\label{fig-psl2}
             Brillouin zones for $PSL(2,\Z)$ in the hyperbolic disk. We
	     have transported the ``usual'' upper half-plane representation
	     using the map $z \mapsto \frac{iz+1}{z+i}$. On the left is
	     shown the sets $B_n(\frac{i}{4})$, which give fundamental
	     domains as in Cor.~\ref{fundamental}.  On the right, $0$ is
	     taken as a basepoint.  Since the origin has a non-trivial
	     stabilizer, the corresponding Brillouin zones give a double
	     cover of the fundamental domains.}
\end{figure}

\begin{cory} 
\label{fundamental}
Let $X$ be $\R^n$, $\S^n$, or $\H^n$, and  let $G$  be a discontinuous group
of  isometries of  $X$. Let $x_0\in X$ be such that its stabilizer $G_{x_0}$
under $G$ is trivial.  Then for every positive integer $n$,  the $n$-th
Brillouin zone  $B_n(x_0)$ is a fundamental domain for the action
of $G$ on points in $X$. Its boundary is the union of
pieces of totally geodesic subspaces and equals the boundary of its interior.
\end{cory}

\remark The first Brillouin zone $B_1(0)$ is the usual Dirichlet 
fundamental domain for the action of $G$.  Furthermore, even when $G_{x_0}$ is not
trivial, $B_n(x_0)$ is a $k$-fold cover of a fundamental domain.

\bigskip
As pointed out in the introduction, the number of geodesics that focus in a
certain point is counted in the lift. So if a given  point $x\in X$ is
intersected  by $n$ mediatrices, it is reached by $n+1$ geodesics of length
$d(0,x)$  emanating from the reference point (the origin). In the next
section, we give more specific examples of this. 
 
Finally, we state a conjecture.
 
\noindent
\begin{conj}
Let $(X,{\tilde g})$ be the universal cover of a $d$-dimensional smooth
Riemannian manifold $(M,g)$ as described in the construction.
For a generic metric $g$ on $M$, no more than $d$ mediatrices
intersect in any given point $y$ of $X$.
\end{conj}
 
This conjecture acquires perhaps even more interest (and certainly
more structure), when one restricts the collection of metrics on
$M$ to conformal ones (\cite{Mas}). A result in this direction for
$M = \R^2/\Z^2$ can be found in \cite{G.A.Jones}.


\section{Focusing in two Riemannian examples}
\label{examples}
\setcounter{figure}{0}
\setcounter{equation}{0}

\bigskip

In this section, we give two examples (one of them new as far as we know)
of focusing. Suppose that at $t=0$ geodesics start emanating in all possible
directions from the a point. At certain times $t_1, t_2, ....$, we may
see geodesics returning to that point. We derive expressions for the number
of geodesics returning at $t_n$ in two cases. First, as introductory example
we will discuss this for $M=\R^2/\Z^2$ (a more complete discussion of this
example can be found in \cite{Pe3}). Second, we will deal with a much more 
unusual example, namely $M=\H^2/\Gamma(k)$, where $\Gamma(k)$ is a
subgroup of $PSL(2,\Z)$ called the principal congruence subgroup of
level $k$ (defined in more detail below).  We note that it seems to be
considerably harder to count geodesics that focus in points other than
our basepoint. 
 
Before continuing, consider the classical problem of counting
$R_g(n)$, the number of solutions in $\Z^2$ of
\bsenn
p^2 + q^2 = n   .
\esenn
Let
\bsenn
n = 2^\alpha \prod_{i=1}^k p_i^{\beta_i} \prod_{j=1}^\ell q_i^{\gamma_i}
\esenn
be the prime decomposition of the number $n$, where $p_i \equiv 1 (\mod\, 4)$
and $q_i \equiv 3 (\mod\, 4)$. The following classical result of Gauss (see
\cite{NZM}) will be very useful.
 
\begin{lem}
$R_g(n)$ is zero whenever $n$ is not an integer, or any of the $\gamma_i$ is
odd.  Otherwise, 
\bsenn
R_g(n) = 4 \prod_{i=1}^k (1+\beta_i)   .
\esenn
\end{lem}

\begin{exam} \rm
\label{flat-torus-example}
Choose an origin in $M=\R^2/\Z^2$ and lift it to the 
origin in $\R^2$. Our discrete set $S$ is then $\Z^2$. Let $\rho_x(t)$ be
the number of geodesics of length $t$ that connect the origin to the point
$x \in M$. 

\begin{prop}
\label{count1}
In the flat torus, the number of geodesics of length $t$ that connect the
the origin to itself is given by 
\bsenn
	\rho_0(t) = \cases{ R_g(t^2), &if $t^2 \in \N$ \cr
				0,    &otherwise.}
\esenn
\end{prop}

\proof Notice that by definition geodesics of length $t$ leaving from
the origin in $\R^2$ reach the points contained in $C_t(0)$. Only if $t^2$
is an integer does this circle intersect points in $\Z^2$. 
\QED
\end{exam}

\begin{exam}\rm\label{gamma2-example}
We now turn to the next example. Recall that $PSL(2,\Z)$ 
can be identified with the group of two by two matrices with integer entries 
and determinant one, and with multiplication by $-1$ as equivalence. 
For each $k$, the group $\Gamma(k)$ is the subgroup of $PSL(2,\Z)$ given by
\bsenn
\Gamma(k) = \set{ \left( \begin{array}{ll} a & b\\ c & d \end{array} \right) 
\in PSL(2,\Z) \st  a \equiv d \equiv 1\, (\mod\, k),\,\, b\equiv c
\equiv 0\, (\mod \,k), }   . 
\esenn
This group has important applications in number theory. The action of
$\Gamma(k)$ on $\H^2$ is given by the M\"obius transformations 
\bsenn
     g(z) = \frac{az+b}{cz+d}   
    \quad\text{where}\quad
    \left( \begin{array}{ll} a & b\\ c & d \end{array} \right) \in \Gamma(k).
\esenn
We point out that for $k=2$, $3$, or $5$, the surface $\H^2 /
\Gamma(k)$ is a sphere with 3, 4, or 12 punctures (see \cite{FK:mod2}).

We will find it more convenient to work in the hyperbolic disk $\D^2$, which
is the universal cover of $\H^2 / \Gamma(k)$.  We shall choose a
representation of $\Gamma(k)$ in the disk so that $i \in \H^2$ corresponds
to the origin.  This will allow us to determine the focusing of the
geodesics which emanate from $i$.  Note that the surface $\H^2 / \Gamma(k)$
has special symmetries with respect to $i$: for example, $i$ is the unique
point fixed by the order $2$ element of $PSL(2,\Z)$.

\begin{figure}[ht]
	\centerline{\psfig{figure=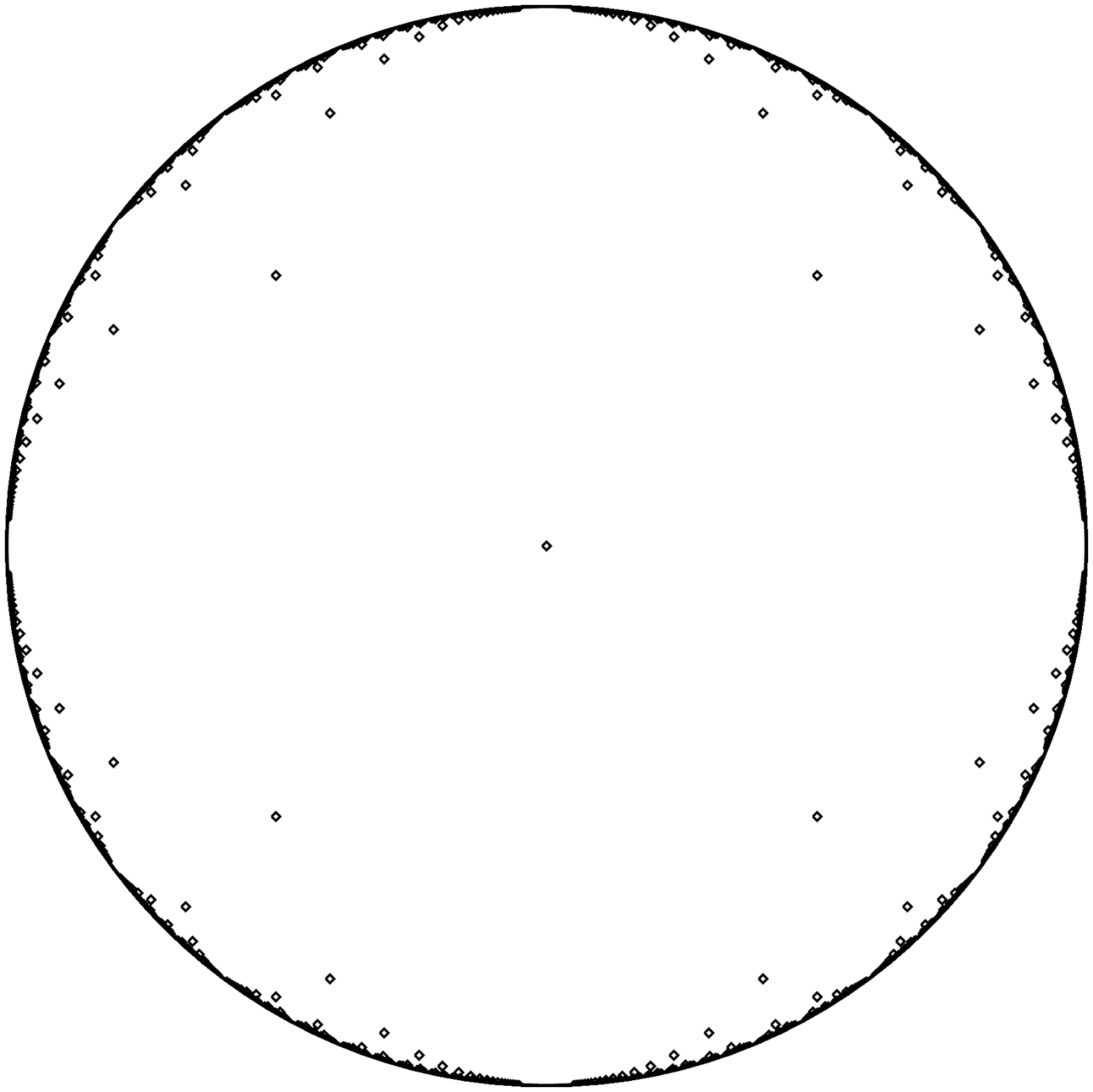,width=.4\hsize} \hfil
		    \psfig{figure=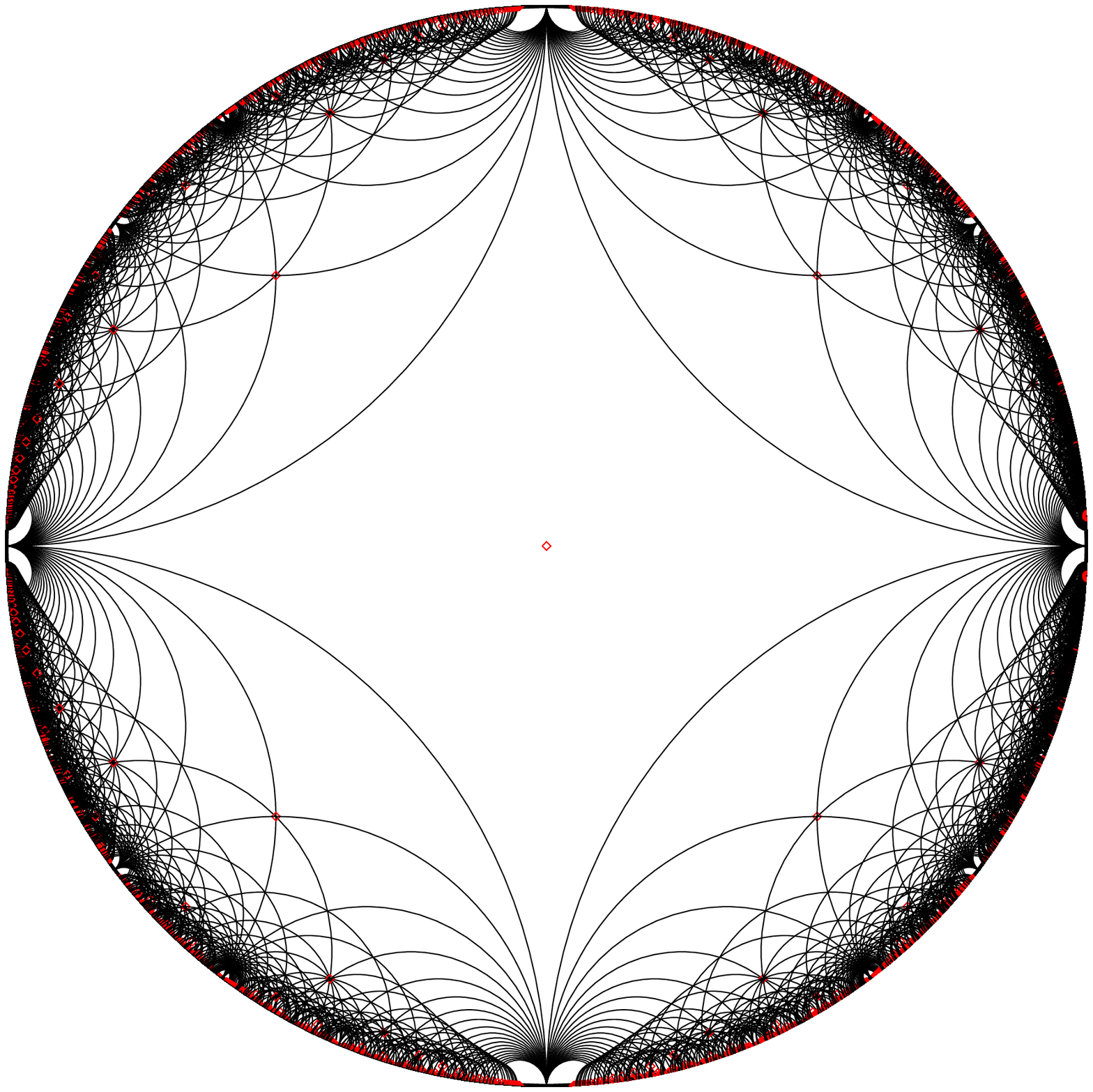,width=.4\hsize}}
	\Caption{The orbit of $i$ under $\Gamma(2)$ transported to the
		 hyperbolic disk, and the corresponding Brillouin zones. 
	         Each zone $B_n$ forms a fundamental domain for a 3-punctured
		 sphere.}
\end{figure}

\begin{lem}
\label{disk_rep}
The action of the fundamental group of the surface $\H^2 / \Gamma(k)$ can be
represented as
\bsenn
\set{ \left( \begin{array}{ll}  r-is & p+iq \cr
                                p-iq & r+is 
             \end{array} \right)
           \st
    \begin{array}{ll}
       r+p \equiv 1 \,(\mod\ k), & r-p \equiv 1 \,(\mod\ k) \cr
       s+q \equiv 0 \,(\mod\ k), & s-q \equiv 0 \,(\mod\ k) 
    \end{array}  
    \text{and}  p^2+q^2+1=r^2+s^2 
},
\esenn
acting on $\D^2$.  We shall denote this particular representation as the
group $\GamDisk{k}$.
\end{lem}


\noindent 
{\bf Proof:} Following the conventions in~\cite{Be}, define
\bsenn
\phi : \D^2 \rightarrow \H^2, \quad \phi(z) = i \frac{z+1}{-z+1}   .
\esenn
Push back the transformation $g \in \Gamma(k)$ from $\H^2$ to $\D^2$ by
$g \rightarrow \phi^{-1}g\phi$ to obtain a representation of $g\in
\Gamma(k)$ as a transformation acting on $\D^2$. The matrix representation
of this transformation is given by:  
\bsenn
A_g = \left( \begin{array}{rr}
 \frac{a+d}{2}+i\frac{b-c}{2} & \frac{a-d}{2} - i\frac{b+c}{2}\\ \ & \  \\
 \frac{a-d}{2}+i\frac{b+c}{2} & \frac{a+d}{2} - i\frac{b-c}{2}
\end{array} \right)   ,
\esenn
where $\det A_g = 1$, since this matrix is conjugate to $g$, whose
determinant is equal to $1$. Let  
\bsenn
\begin{array}{rcr}
	p=(a-d)/2  & \qquad & q= -(b+c)/2 \\
	r=(a+d)/2 & \qquad & s=-(b-c)/2
\end{array}
\esenn
and $A_g$ now written as
\bsenn
A_g = \left( \begin{array}{lr}
	r-is & p+iq  \\
	p-iq & r+is \end{array} \right)   .
\esenn
Here the numbers $p,q,r,s$ are in $\Z$ and must satisfy the following
congruence conditions:
\bsenn
    \begin{array}{ll}
       r+p \equiv 1 \,(\mod\ k), & r-p \equiv 1 \,(\mod\ k) \cr
       s+q \equiv 0 \,(\mod\ k), & s-q \equiv 0 \,(\mod\ k) 
    \end{array}  
\esenn
Since the determinant of $A_g$ is equal to $1$, we must also have 
	$$p^2+q^2+1=r^2+s^2.$$
\QED

We need another auxiliary result before we state the main result of this
section. 

\begin{lem}
\label{auxiliary}
Let $(p,q)$ and $(r,s)$ be two points in $\Z^2$ such that the integers
$A=p^2 + q^2$ and $B=r^2+s^2$ are relatively prime, and let $\varphi$ be a
rotation. 
Now $\varphi(p,q)=(p',q')$ and $\varphi(r,s)=(r',s')$ are in $\Z^2$ if and
only if $\varphi$ is a rotation by a multiple of $\pi/2$.
\end{lem}

\proof 
Let $c$ be the cosine of the angle of rotation.
We have 
\bsenn
c= \frac{p'p+q'q}{A}=\frac{r'r+s's}{B}   .
\esenn
Thus if $p'p+q'q$ and $r'r+s's$ are not both equal to zero, 
\bsenn
\frac{p'p+q'q}{r'r+s's} = \frac{A}{B}   .
\esenn
Because $A$ and $B$ are relatively prime and surely $|p'p+q'q|$ is less than
or equal to $A$, and similarly for $B$, we have that 
\bsenn
p'p+q'q= \pm A \logand r'r+s's = \pm B   .
\esenn
This implies the result.
\QED

Now we define a counter just as before.
Choose a lift of $M = \D^2 / \GamDisk{k}$ so that  $0\in M$ lifts to
$0\in\D^2$. Let $\gamma_x(t)$ be the number
of geodesics of length $t$ that connect the origin to the point $x \in M$.

\begin{figure}[ht]
	\centerline{\psfig{figure=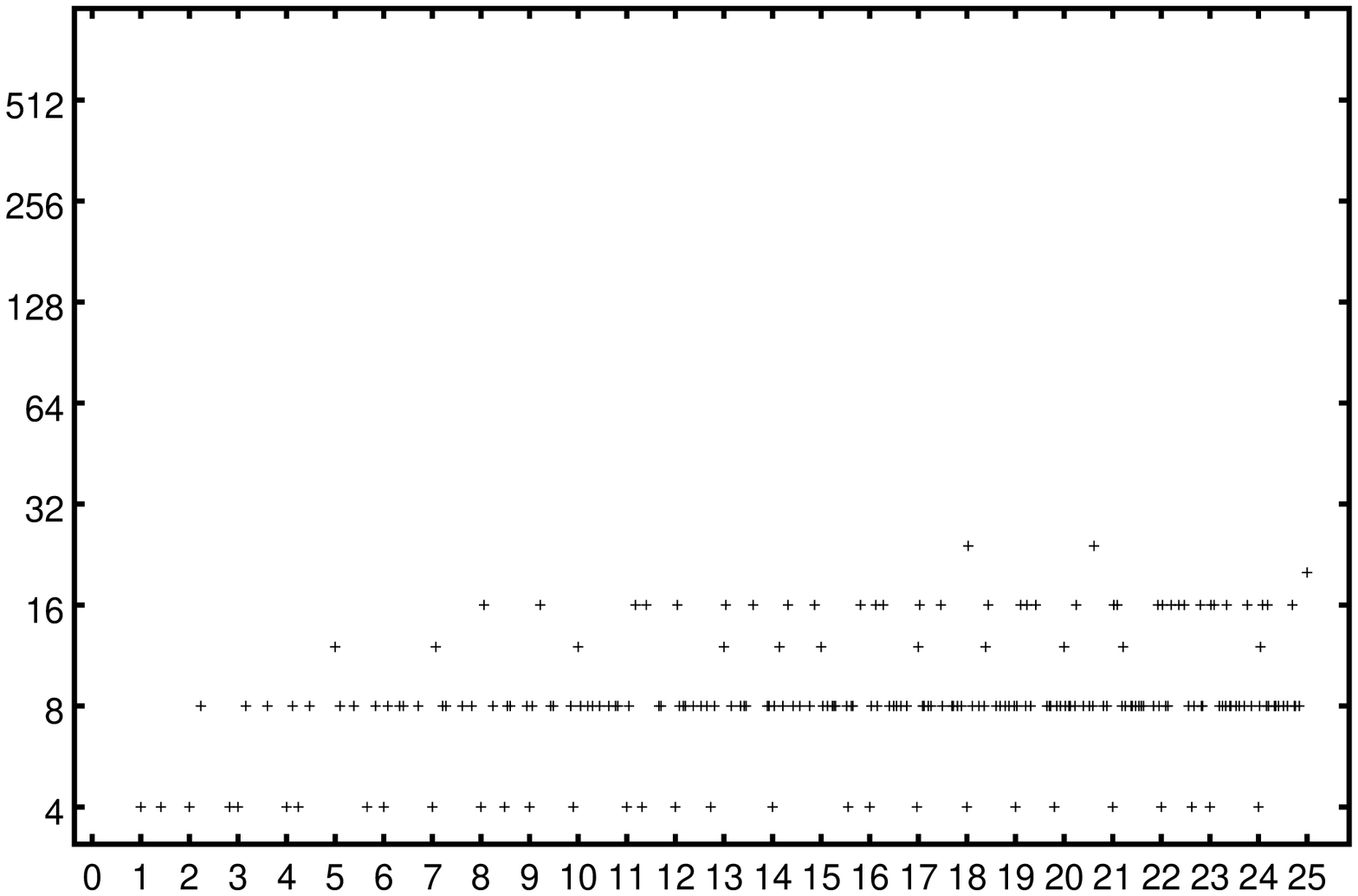,height=.4\hsize} \quad\quad
                    \psfig{figure=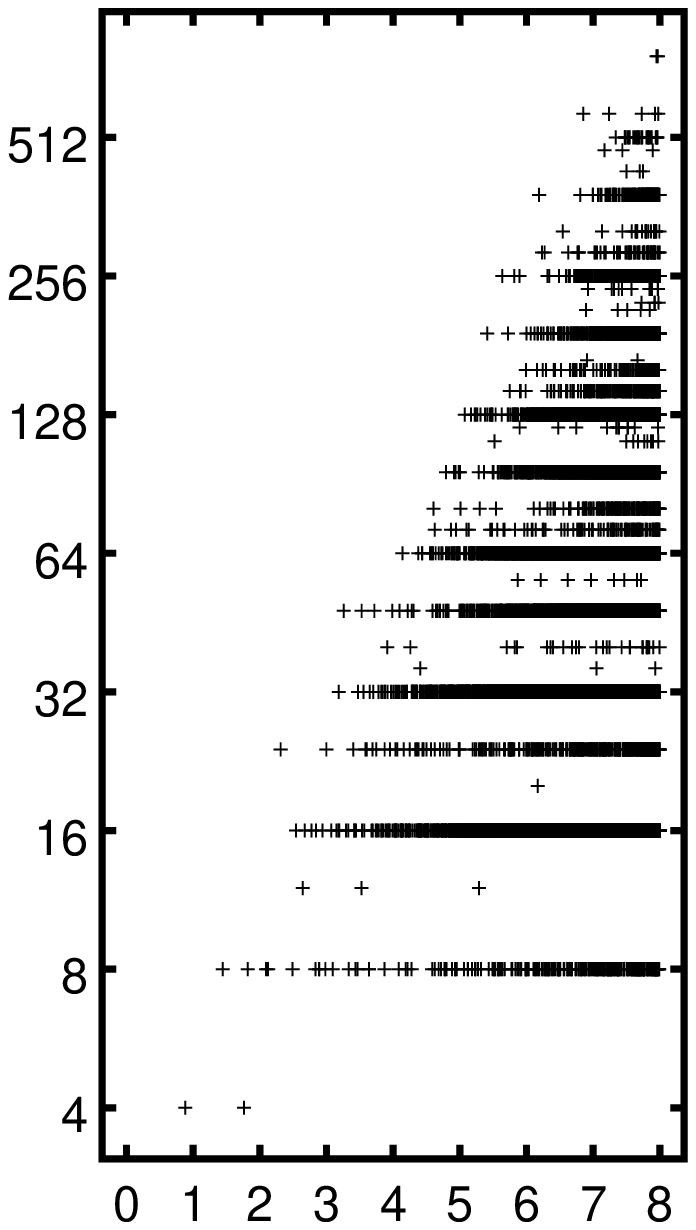,height=.4\hsize}}
	\Caption{The non-zero values of $\rho_0(t)$ for $t\le 25$ (left) and
                    $\gamma_0(t)$ for $t \le 8$ (right), which count how
                    many geodesics of length $t$ connect the origin to
                    itself in the $\R^2 / \Z^2$ and $\D^2 / \GamDisk{2}$,
                    respectively.} 
\end{figure}

\begin{theo}
\label{gamma_k}
In the surface $\H^2 / \Gamma(k)$, the number of geodesics which connect the
point $i\in\H^2$ to itself is given by
\bsenn
 \begin{array}{lcr}
  \frac{1}{4} R_g(\cosh^2 t -1 ) R_g(\cosh^2 t)     & & \text{for} k=2 \\[3pt]
  \frac{1}{4} R_g \left( \frac{\cosh^2 t -1}{9} \right) R_g(\cosh^2 t)
						     & & \text{for} k=3 \\[5pt]
  \frac{1}{4} R_g \left( \frac{\cosh^2 t -1}{25} \right) R_g(\cosh^2 t)
						     & & \text{for} k=5 \\
 \end{array}
\esenn
Note that in all cases, the number is nonzero only if $\cosh^2 t \in \N$.
\end{theo}

\proof
We shall work in the disk, rather than in $\H^2$.
Let $S$ be the orbit of $0\in \D^2$ under $\GamDisk{k}$.
Then the number of such geodesics is exactly the number of distinct points
of $S$ which lie on on the circle  $C_t(0)$ of radius $t$ and centered at the
origin.

If $x\in S$, then by Lemma~\ref{disk_rep}, it is of the form
$\frac{p+iq}{r+is}$ with $p$, $q$, $r$, $s$ integers satisfying 
$p^2 + q^2 = r^2 + s^2 -1$.  Let $n$ be their common value, that is,
\bsenn
n = p^2+q^2 = r^2+s^2-1.
\esenn
We will first count the number of 4-tuples $(p,q,r,s)$ that solve this 
equation, momentarily ignoring the congruence conditions.

\goodbreak

Note that the point $x$ has Euclidean distance to the origin given by 
\bsenn
|x|_e^2 = \frac{p^2+q^2}{r^2+s^2}  = \frac{n}{n+1} .
\esenn
The hyperbolic length of the geodesic connecting the origin to this point is
\bsenn
t = \atanh(|x|_e)   .
\esenn
This gives that $\gamma_0(t)$ is only non-zero when
$t=\atanh \sqrt{\frac{n}{n+1}}$, or, equivalently, 
$n=\cosh^2 t-1$.

To count the number of intersections of $C_t(0)$ with $S$ for these
values of $t$, observe that we can use Gauss' result to 
count the number of pairs $(p,q)$ such that $p^2+q^2=n$. This number
is given by $R_g(n)$. For each such pair $(p,q)$, we have a number of
choices to form 
\bsenn
x=\frac{p+iq}{r+is}
\esenn
By the above, this number is equal to $R_g(n+1)$.  
Thus, $\gamma_0(t)$ is at most $R_g(n) R_g(n+1)$.  


However, we have over-counted: some of our choices for $p,q,r,s$ represent
the same point $x \in S$, and some of them may not satisfy the congruence
conditions, which we have so far ignored. We will first account for the
multiple representations, and then account for the congruence relations.

Let $p,q,r,s \in \Z$ be as above, giving a point $x=\frac{p+iq}{r+is}$ at
distance $t=\atanh\sqrt{\frac{n}{n+1}}$ from the origin. If we multiply the
numerator and denominator of $x$ by $e^{i\theta}$, then $x$ will remain
unchanged.  Because of the requirement that $p^2+q^2=r^2+s^2-1=n$, this is
the only invariant, and by Lemma~\ref{auxiliary}, $\theta$ must be a
multiple of $\frac{\pi}{2}$ for the numerator and denominator to remain
Gaussian integers.  We see that in our counting, we have represented our
point $x$ in 4 different ways:
\bsenn
x = \frac{ p+iq}{ r+is} = \frac{-q+ip}{-s+ir} = 
    \frac{-p-iq}{-r-is} = \frac{ q-ip}{ s-ir}
\esenn
meaning we have over-counted by a factor of at least 4.  

Now we account for the congruence conditions.  

First, consider the case $k=2$.
Note that $q+s \equiv 0 \,(\mod\ 2)$ if and only if $p+r \equiv 1 \,(\mod\
2)$, because $p^2 + q^2 + 1 = r^2 + s^2$, so we need only check this one
condition. If the representation $\frac{p+iq}{r+is}$ fails to
satisfy our parity condition, then $q+s \equiv 1 \,(\mod\,2)$ and consequently
$p+r \equiv 0 \,(\mod\,2)$. This means that the representation
$\frac{-q+ip}{-s+ir}$ of this same point does satisfy the parity conditions,
giving exactly $\frac{1}{4}R_g(\cosh^2 t -1 ) R_g(\cosh^2 t)$ distinct
points of $S$ at distance $t$ from the origin. 

\medskip
If $k=3$, then since $k$ is odd, the congruence conditions on $p,q,r,$ and
$s$ imply that  
\bsenn
  r \equiv 1 (\mod\ 3)    \quad\logand\quad 
  p \equiv q \equiv s \equiv 0 (\mod\ 3)  .
\esenn
Note that the equation
\bsenn
  p^2 + q^2 = n   \quad\logand\quad    p \equiv q \equiv 0 \,(\mod\ 3)
\esenn
will be satisfied exactly $R_g(n/3^2)$ times. (Recall that if $n/9$ is not an
integer, $R_g(n/9) = 0$.)

\goodbreak
For fixed $n$, let $(p,q)$ be any one of the solutions.  We need to decide
how many solutions the equation 
\bsenn
  r^2+s^2=n+1 \quad\text{with}\quad
  r\equiv 1 \,(\mod\ 3) \quad\logand\quad s\equiv0 \,(\mod\ 3)
\esenn
admits. The solution of the first equation implies that 3 divides $n$. 
Thus $r^2+s^2 \equiv 1 \,(\mod\ 3)$.
Consequently, we have 4 choices $\mod 3$ for the pair $(r,s)$, namely
$(0,1)$, $(1,0)$, $(0,2)$, and $(2,0)$.

Let $(p,q,r,s) \in \Z^2 \times \Z^2$ be any solution to 
$n=p^2+q^2 = r^2 + s^2 -1$ with $p \equiv q \equiv 0 (\mod\ 3)$.
For each choice of $(p,q)$, we have exactly $R_g(n+1)$ choices of $(r,s)$.
Now let $R$ denote the product of the rotations by $\pi/2$ on each of the
components of $\Z^2 \times \Z^2$.  Using Lemma~\ref{auxiliary}, we see that
all such solutions can be obtained from just one by applying $R$ repeatedly. 
It is easy to check that each quadruple of solutions thus constructed runs
exactly once through the above list. Since precisely one out of the four
associated solutions is compatible with the conditions,
the total number of solutions is exactly:
\bsenn
{1 \over 4} R_g \left({n \over 9} \right) R_g(n+1)  .
\esenn
Using the relationship between the Euclidean distance and the Poincar\'e
length as before gives the result.

\medskip
If $k=5$, the proof for $k=3$ can be literally transcribed to obtain the
result. 
\QED

\remark
Note that the above results do not hold if $k$ is not one of the cases
mentioned.  The primary difficulty is that for prime $k \ge 7$, there are
solutions which are not related by applying the rotation $R$.  However, the
argument does give 
 $\frac{1}{4} R_g \left( (\cosh^2 t -1)/k^2 \right) R_g(\cosh^2 t) $
as an upper bound for $\H^2/\Gamma(k)$ when $k$ is an odd prime.
Note that the surface $\H^2/\Gamma(k)$ is of genus $0$ if and only if 
$k \le 5$ (see \cite{FK:mod2}).

\end{exam}

\section{Non-Riemannian examples}
\label{others}
\setcounter{figure}{0}
\setcounter{equation}{0}

\bigskip

The present context is certainly not restricted to Riemannian metrics. As
an indicator of this we now discuss a different set of examples.

Let $k$ be a positive number greater than one. Equip $\R^2$ with the 
distance function 
\bsenn
 \parallel x-y \parallel = \left( |x_1-y_1|^k + |x_2-y_2|^k \right)^{1/k}   
\esenn
and let the discrete set $S$ be given by $\Z^2$.
For $k$ not equal to 2, this is not a Riemannian metric, yet all conclusions
of section \ref{Main Results} hold.   In particular, each Brillouin zone
forms a fundamental domain.  Note that determining the zones by inspecting
the picture requires close attention!

\begin{figure}[ht]
 \centerline{\psfig{figure=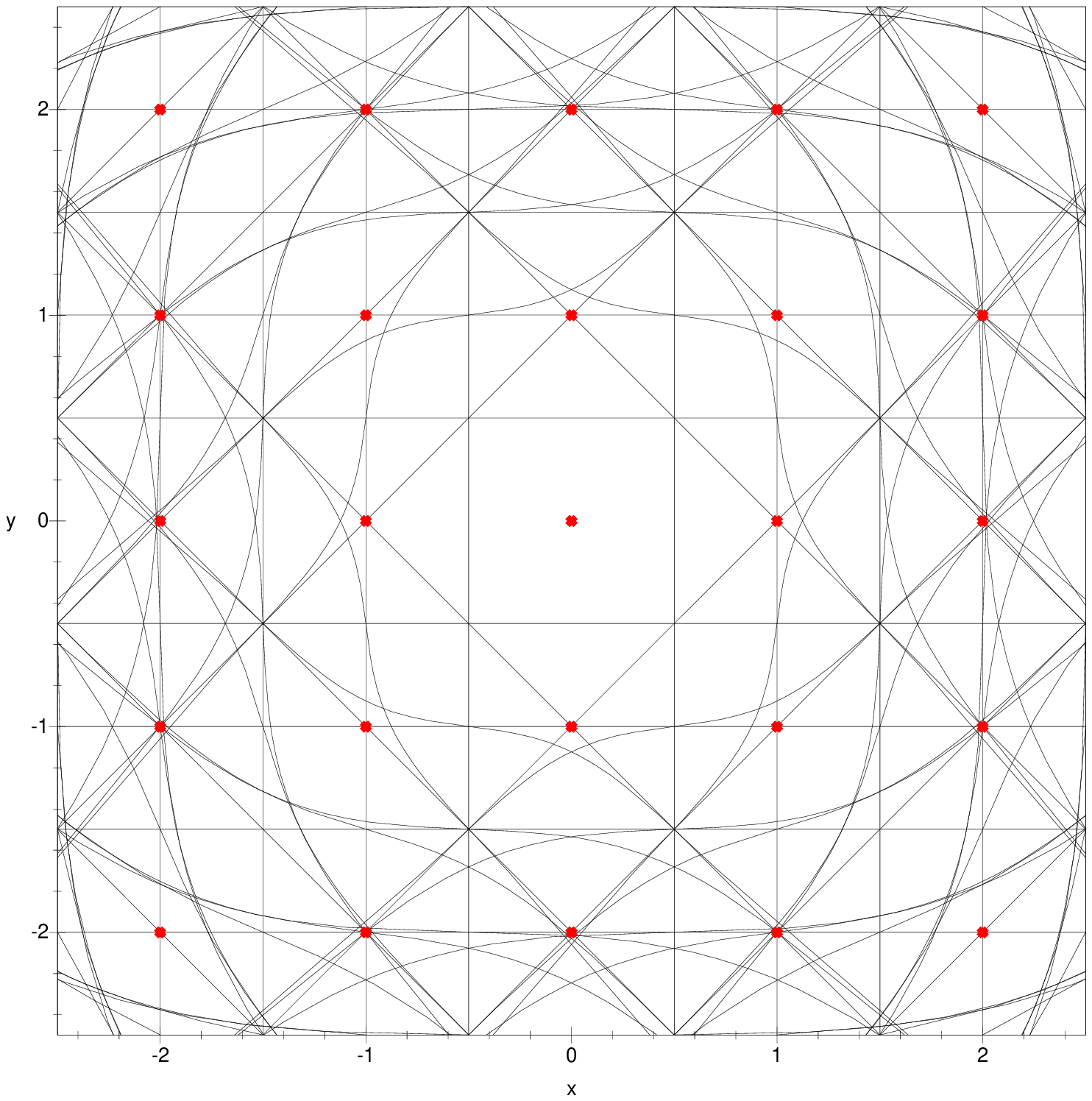,width=.45\hsize} \hfil
             \psfig{figure=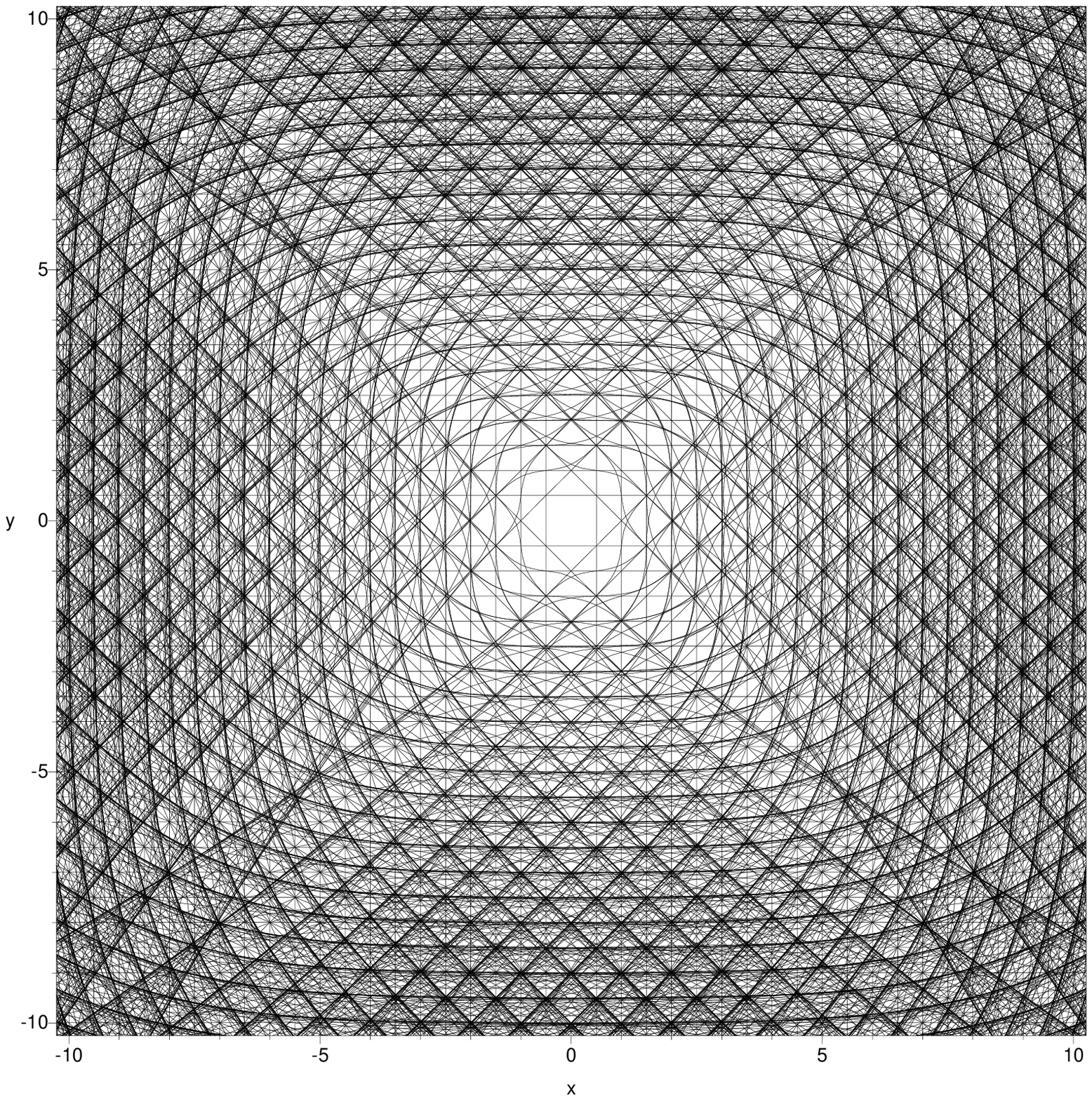,width=.45\hsize}}
 \Caption{\label{lattice4-fig} 
     Brillouin zones for the lattice $\Z^2$ in $\R^2$ with the 
     metric $ \left(|x_1-y_1|^4 + |x_2-y_2|^4 \right)^{1/4}$.
     See also Figure~\ref{fig-man3} and Example~\ref{ex-manhattan}, which
     deal with the case $k=1$, the ``Manhattan metric''.
}
\end{figure}

Now the problem of determining $C_t(0)\intersect S$ for any
given $t$ is unsolved for general $k$. In fact, even for
certain integer values of $k$ greater than 2, it is not known
whether $C_t(0)\intersect S$ ever contains at least two points
that are not related by the symmetries of the problem.
For $k=4$, the smallest $t$ for which $C_t(0)\intersect S$
has at least two (unrelated) solutions is given by
\bsenn
   t^4=133^4+134^4=158^4+59^4     .
\esenn
However, for $k \ge 5$, it unknown whether this can happen at all
(see \cite{SW}).

There are some things that can be said, however. In the
situation where $k\in\set{3, 4, 5, \ldots}$, the mediatrices intersect
the coordinate axes only in irrational points or in multiples
of $1/2$. For if $x=(p/q,0)$ is a point of a mediatrix
$L_{(0,0),(a_1,a_2)}$, we have
\bsenn
|p|^k = |p - q a_1|^k +|a_2 q|^k \quad (p\ne 0, q\ne 0) \quad .
\esenn
By Fermat's Last Theorem, this has no solution unless either
$p = q a_1$ or $a_2= 0$.  In the first case, ${p \over q} = \pm a_2$, which can
only occur if the lattice point is of the form $(a_2, \pm a_2)$.  If $a_2=0$,
then ${p \over q}={a_1 \over 2}$. 
In particular, there is no nontrivial focusing along the axes.

\medskip
To compute Figure~\ref{lattice4-fig}, we took advantage of the smoothness of
the metric.  Not all metrics are sufficiently smooth for this procedure to
work.  Even for Riemannian metrics, in general the distance function is
only Lipschitz, which will not be sufficiently smooth.

For each $a=(a_1,a_2) \in \Z^2$,  define a Hamiltonian:
\bsenn
H_a(x) = \parallel x-a \parallel - \parallel x \parallel   .
\esenn
The mediatrix $L_{0a}$ corresponds to the level set $H_a(x) = 0$.
Because $H_a(x)$ is smooth, we have uniqueness of solutions to
Hamilton's equations. In the current situation, where the dimension is two, 
the level set consists of one orbit. Thus, one can produce the mediatrix
by numerically tracing the zero energy orbits of the above Hamiltonian.

\medskip

As mentioned above, for a general Riemannian metric, the distance function
is only Lipschitz. This means we have no guarantee that the solutions of the
above differential equation are unique. Indeed, there are examples of
multiply connected Riemannian manifolds with self-intersecting mediatrices,
as will be shown in a forthcoming work. 

\bigskip\goodbreak

\end{document}

%% file: inputlatex.tex
\newtheorem {theo} {\bf Theorem} [section]
\newtheorem {prop} [theo] {\bf Proposition}
\newtheorem {cory} [theo] {\bf Corollary}
\newtheorem {lem} [theo] {\bf Lemma}
\newtheorem {defn} [theo] {\bf Definition}
\newtheorem {conj} [theo] {\bf Conjecture}
\newtheorem {exam} [theo] {\bf Example}
\newtheorem {rem} [theo] {\bf Remark}
\newcommand{\QED}{\hfill \CaixaPreta \vspace{6mm}}

\newenvironment{remark}{\begin{rem}\begin{rm}}{\end{rm}
\nopagebreak\hfill$\CaixaBranca$\end{rem}\par\vspace{\partopsep}}

\newcommand{\qed}{\nopagebreak\hfill{\vrule width6pt height6pt depth0pt}}
\newenvironment{proof}[1]{\vspace{0.4 cm}\noindent{\em #1.}}{\qed\par
\vspace{\topsep}\vspace{\partopsep}}

\newcommand{\be}{\begin{eqnarray}}
\newcommand{\ee}{\end{eqnarray}}
\newcommand{\benn}{\begin{eqnarray*}}
\newcommand{\eenn}{\end{eqnarray*}}
\newcommand{\bse}{\begin{equation}}
\newcommand{\ese}{\end{equation}}
\newcommand{\bsenn}{\begin{displaymath}}
\newcommand{\esenn}{\end{displaymath}}
\newcommand{\logand}{\;\;{\rm and }\;\;}

\newcommand{\with}{\;\;{\rm with }\;\;}

\newcommand{\atanh}{\;{\rm arctanh }\,}

\newcommand{\mod}{\bmod}
\renewcommand{\H}{\mbox{$\Bbb H$}} 	
\newcommand{\N}{\mbox{$\Bbb N$}} 	
\newcommand{\R}{\mbox{$\Bbb R$}} 	
\newcommand{\D}{\mbox{$\Bbb D$}} 	
\newcommand{\Z}{\mbox{$\Bbb Z$}} 	
\renewcommand{\S}{\mbox{$\Bbb S$}}  
\newcommand{\CaixaBranca}{$\square$}


%% file: imsmark.tex
\def\SBIMSMark#1#2#3{
 \font\SBF=cmss10 at 10 true pt
 \font\SBI=cmssi10 at 10 true pt
 \setbox0=\hbox{\SBF Stony Brook IMS Preprint \##1}
 \setbox2=\hbox to \wd0{\hfil \SBI #2}
 \setbox4=\hbox to \wd0{\hfil \SBI #3}
 \setbox6=\hbox to \wd0{\hss
             \vbox{\hsize=\wd0 \parskip=0pt \baselineskip=10 true pt
                   \copy0 \break%
                   \copy2 \break%
                   \copy4 \break}}
 \dimen0=\ht6   \advance\dimen0 by \vsize \advance\dimen0 by 8 true pt
                \advance\dimen0 by -\pagetotal
 \dimen2=\hsize \advance\dimen2 by .25 true in
 \ht6=0pt \dp6=0pt
 \setbox8=\vbox to \dimen0{\vfill \hbox to \dimen2{\hss \copy6}}
 \ht8=0pt \dp8=0pt \wd8=0pt
 \copy8
 \message{*** Stony Brook IMS Preprint #1, #2 ***}
}